%% file: main.tex
\documentclass{amsart}
\usepackage{graphicx,amsmath,amsthm,amssymb,amsfonts,mathrsfs,quiver,url} 
\usepackage{tikz-cd}

\title[Equivariant geometry of quadrics blown-up along a twisted quartic]{Equivariant geometry of quadric threefolds blown up along a twisted quartic}
\author{Joseph Malbon, Antoine Pinardin}
\date{October 2026}

\include{antoinescommands.tex}

\newcommand{\Gm}{\mathbb{G}_m}

\newcommand{\wt}[1]{\widetilde{#1}}
\newcommand{\oldphi}{\mathchar"011E}
\DeclareMathOperator{\Bs}{Bs}
\numberwithin{equation}{section}
\newcommand{\C}{\mathbb{C}}

\begin{document}
\maketitle
\begin{abstract}
    Let $X$ be the blow-up of a smooth quadric threefold along a twisted quartic, and let $G\subset\aut(X)$ be a finite subgroup such that $\pic^G(X)$ has rank one. We prove that $X$ is $G$-birationally superrigid if and only if $G$ is isomorphic to $\ac\times\dz_2$. The action of $G$ on $X$ is linearisable if and only if it is projectively linearisable, and this holds precisely when $G$ contains no subgroup isomorphic to $\dz_2^3$.
\end{abstract}
\section{Introduction}
We work over the field of complex numbers. Let $X$ be a smooth Fano threefold in the family \textnumero 2.21, namely the blow-up $\pi\colon X\rightarrow Q$ of a smooth quadric threefold $Q\subset\p^4$ along a smooth rational normal quartic curve $C_4$, also called a \textit{twisted quartic}. For a finite subgroup $G\subset\operatorname{Aut}(X)$, questions of $G$-birational rigidity and linearisability arise when we ask which $G$-Mori fibre spaces realise its conjugacy class in the Cremona group of rank three and, in particular, whether this class contains a linear action on $\p^3$. We use the following description of $\aut(X)$.

\begin{theo}[\cite{malbon2024automorphisms}]
\label{theo:aut-pres} Let $X$ be a smooth Fano threefold in the family \textnumero 2.21. There exists an involution $\tau\in\aut(X)$ such that $\aut(X)\cong\aut(Q,C_4)\times\langle\tau\rangle$.
\end{theo}

\noindent We identify $\aut(Q,C_4)$ with its lift to $\aut(X)$ whenever no confusion can arise. Hence, there is a commutative diagram

$$
\begin{tikzcd} & X & \\ {Q} && {Q} \arrow["\pi"', from=1-2, to=2-1] \arrow["{\pi'}", from=1-2, to=2-3] \arrow["\sigma"', dashed, from=2-1, to=2-3] \end{tikzcd}
$$
\noindent where $\pi'=\pi\circ\tau$ and $\sigma=\pi'\circ\pi^{-1}$. We denote the exceptional divisors of $\pi$ and $\pi'$ by $E$ and $E'$, respectively. The group $\aut(Q,C_4)$ preserves both extremal contractions of $X$, while $\tau$ exchanges them. Thus, for a finite subgroup $G\subset\aut(X)$, the invariant Picard group $\pic^G(X)$ has rank one precisely when some element of $G$ exchanges the two contractions.

Recall that a $G$-Fano variety $X$ is called \emph{$G$-birationally superrigid} if every $G$-birational map $X\dashrightarrow Y$ to a $G$-Mori fibre space $Y\to S$ has $S$ a point and is biregular; it is \emph{$G$-birationally rigid} if the same holds after precomposition with a $G$-birational self-map of $X$. It is \emph{$G$-birationally solid} if there is no such map with $\dim S>0$. We abbreviate these terms to \emph{$G$-superrigid}, \emph{$G$-rigid}, and \emph{$G$-solid}.

The action on our threefold $X$ is \emph{linearisable} if $X$ is $G$-birational to $\mathbb P(V)$ for a four-dimensional linear representation $V$ of $G$; it is \emph{projectively linearisable} if $X$ is $G$-birational to $\mathbb P^3$ endowed with a biregular $G$-action. The main results of this work are the following.

\begin{theo}\label{theo:rigid}
Let $X$ be a Fano threefold in the family \textnumero2.21, and $G$ be a finite subgroup of $\aut(X)$ such that $\pic^G(X)\cong\dz$. Then the following conditions are equivalent.
\begin{enumerate}
\item $G$ is isomorphic to $\ac\times\dz_2$,
\item $X$ is $G$-birationally rigid,
\item $X$ is $G$-birationally superrigid,
\item $X$ is $G$-birationally solid.
\end{enumerate}
\end{theo}

\begin{theo}\label{theo:lin}
Let $X$ be a Fano threefold in the family \textnumero2.21, and $G$ be a finite subgroup of $\aut(X)$ such that $\pic^G(X)\cong\dz$. Then the following conditions are equivalent.
\begin{enumerate}
\item the action of $G$ on $X$ is linearisable,
\item the action of $G$ on $X$ is projectively linearisable,
\item $G$ does not contain any subgroup isomorphic to $\dz_2^3$.
\end{enumerate}
\end{theo}

\section{The automorphism group of $X$}

\subsection{The stabiliser of $C_4$ in $\aut(Q)$}

After a projective change of coordinates in $\mathbb P^4$, the curve $C_4$ is the rational normal quartic curve of $\mathbb P^4$, namely the image of $\mathbb P^1$ via the Veronese embedding
$$
\mathcal V\colon [u:v]\longmapsto[u^4:u^3v:u^2v^2:uv^3:v^4].
$$
We use these coordinates on $Q$ throughout, unless otherwise stated. The following polynomials generate $\mathcal I_{C_4}(2)$.
\begin{align*}
f_0 &= x_3^2 - x_2x_4, \\ f_1 &= x_2x_3 - x_1x_4, \\ f_2 &= x_2^2 - x_0x_4, \\ f_3 &= x_1x_2 - x_0x_3, \\ f_4 &= x_1^2 - x_0x_2, \\ f_5 &= 3x_2^2 - 4x_1x_3 + x_0x_4.
\end{align*}

Restriction to $C_4\cong\pl$ identifies $\aut(Q,C_4)$ with a subgroup of $\pgl_2(\dc)$. We give generators in these coordinates. For the classification in the cases with infinite automorphism group, see also \cite{cheltsov2018fano}.

\begin{theo}[\cite{malbon2024automorphisms}]\label{theo:autQC}
Up to an automorphism of $\mathbb P^4$, one of the following holds.
\begin{enumerate}\label{joe:actions}
\item $Q=\{f_5=0\}$, $\aut(Q,C_4)$ is isomorphic to $\pgl_2(\dc)$.\label{joe1}
\item $Q=\{\lambda f_2+f_5=0\}$, for some $\lambda\in\dc\setminus\{0,1,-3\}$, $\aut(Q,C_4)$ is isomorphic to $\dg_m\rtimes\dz_2$, generated by

$$
\dg_m=\left\{\begin{pmatrix}1&0\\0&t\end{pmatrix}\,\middle|\, t\in\dg_m\right\} \qquad\text{and}\qquad \dz_2=\left\langle \begin{pmatrix}0&1\\1&0\end{pmatrix}\right\rangle.
$$
\label{joe2}
\item $Q=\{f_0+f_5=0\}$, $\aut(Q,C_4)$ is isomorphic to $\dg_a\rtimes\dz_2$, generated by

$$
\dg_a=\left\{\begin{pmatrix}1&t\\0&1\end{pmatrix}\,\middle|\, t\in\dg_a\right\} \qquad\text{and}\qquad \dz_2=\left\langle \begin{pmatrix}1&0\\0&-1\end{pmatrix}\right\rangle.
$$
\label{joe3}
\item $Q=\{\mu(f_0+f_4)+\lambda f_2+f_5=0\}$, for some $\lambda\in\dc\setminus\{1,-3\}$ and $\mu\in\dc\setminus\{0,2,-2\}$ such that $\mu^2\ne-\lambda^2-2\lambda+3$ and $\mu^2\ne4\lambda^2$, $\aut(Q,C_4)$ is isomorphic to $\dz_2^2$, generated by
$$g_1=\begin{pmatrix}1&0\\0&-1\end{pmatrix},\qquad g_2=\begin{pmatrix}0&1\\1&0\end{pmatrix}.$$
\label{joe4}

\item $Q=\{f_0+\lambda f_2+f_5=0\}$, for some $\lambda\in\dc\setminus\{0,1,-3\}$, $\aut(Q,C_4)$ is isomorphic to $\dz_2$, generated by

$$
g_1=\begin{pmatrix}1&0\\0&-1\end{pmatrix}.
$$
\label{joe5}
\item $Q=\{f_1+f_5=0\}$, $\aut(Q,C_4)$ is trivial.\label{joe6}
\end{enumerate}
\end{theo}

\subsection{Finite subgroups of $\aut(X)$}

For every finite subgroup $G\subset\aut(X)$, put $H=G\cap \aut(Q,C_4)$. Let $\zeta_n=\exp(2\pi i/n)$, and consider the following elements of $\pgl_2(\dc)$:
$$
R_n=\begin{pmatrix}1&0\\0&\zeta_n\end{pmatrix}, \qquad B=\begin{pmatrix}0&1\\1&0\end{pmatrix},
$$
$$
C=\begin{pmatrix}i&-i\\1&1\end{pmatrix}, \qquad D=\begin{pmatrix}1&-i\\i&-1\end{pmatrix},
$$
$$
M=\begin{pmatrix}\zeta_5&0\\0&1\end{pmatrix}, \qquad N=\begin{pmatrix} 1&1-\zeta_5-\zeta_5^{-1}\\ 1&-1 \end{pmatrix}.
$$
Whenever these transformations preserve $Q$, we use the same notation for their lifts to $X$. They satisfy
$$
\langle R_2,B,C\rangle\cong\mathfrak A_4, \qquad \langle R_2,B,C,D\rangle\cong\mathfrak S_4, \qquad \langle M,N\rangle\cong\mathfrak A_5.
$$
We denote by $D_n$ the dihedral group of order $2n$, with the conventions $D_1\cong\dz_2$, $D_2\cong\dz_2^2$, and $\dz_1=\{\mathrm{id}\}$.

\begin{coro}\label{coro:finite-subgroups-autX}
Every finite subgroup $G\subset\aut(X)$ is conjugate to exactly one of the representatives specified below.

\begin{enumerate}
\item In cases~\ref{joe1} and~\ref{joe2}, define the following groups for $n\geq 1$:

$$
\begin{array}{c|c|c|c}
G & \text{Isomorphism type} & H=G\cap \aut(Q,C_4) & \operatorname{rk}\pic^G(X) \\ \hline
G_{n,1}=\langle R_n\rangle & \dz_n & \langle R_n\rangle\cong\dz_n & 2 \\[2pt]
G_{n,2}=\langle R_n,B\rangle & D_n & \langle R_n,B\rangle\cong\dih_n & 2 \\[2pt]
G_{n,3}=\langle R_n,\tau\rangle & \dz_n\times\dz_2 & \langle R_n\rangle\cong\dz_n & 1 \\[2pt]
G_{n,4}=\langle R_{2n}\tau\rangle & \dz_{2n} & \langle R_n\rangle\cong\dz_n & 1 \\[2pt]
G_{n,5}=\langle R_n,B,\tau\rangle & D_n\times\dz_2 & \langle R_n,B\rangle\cong\dih_n & 1 \\[2pt]
G_{n,6}=\langle R_n,B\tau\rangle & D_n & \langle R_n\rangle\cong\dz_n & 1 \\[2pt]
G_{n,7}=\langle R_{2n}\tau,B\rangle & D_{2n} & \langle R_n,B\rangle\cong\dih_n & 1
\end{array}
$$
In case~\ref{joe2}, all groups in this table, with $n\geq 1$, form an irredundant list of representatives. In case~\ref{joe1}, retain $G_{n,1},G_{n,3},G_{n,4}$ for $n\geq 1$ and $G_{n,2},G_{n,5},G_{n,6},G_{n,7}$ for $n\geq 2$. There are additionally

$$
\langle R_2,B,C\rangle,\qquad \langle R_2,B,C,D\rangle,\qquad \langle M,N\rangle,
$$
isomorphic respectively to $\mathfrak A_4$, $\mathfrak S_4$, and $\mathfrak A_5$; the three groups obtained by adjoining $\tau$ to these generating sets; and the group

$$
\mathfrak S^{\mathrm{tw}} =\langle R_2,B,C,D\tau\rangle \cong\mathfrak S_4.
$$
For the last group, we have

$$
H=\langle R_2,B,C\rangle\cong\mathfrak A_4.
$$
\item In cases~\ref{joe3} and~\ref{joe5}, the groups are

$$
\{\mathrm{id}\},\qquad \langle R_2\rangle,\qquad \langle\tau\rangle,\qquad \langle R_2\tau\rangle,\qquad \langle R_2,\tau\rangle,
$$
of isomorphism types $\{\mathrm{id}\},\dz_2,\dz_2,\dz_2,\dz_2^2$, respectively.
\item In case~\ref{joe4}, the groups are precisely

$$
G_U=\left\{R_2^aB^b\tau^c\,\middle|\,(a,b,c)\in U\right\}\cong\dz_2^{\dim U},
$$
where $U$ ranges over the vector subspaces of $\mathbb Z_2^3$.

\item In case~\ref{joe6}, the groups are $\{\mathrm{id}\}$ and $\langle\tau\rangle\cong\dz_2$.
\end{enumerate}
Moreover,
$$
\operatorname{rk}\pic(X)^G=1 \quad\Longleftrightarrow\quad G\not\subset \aut(Q,C_4).
$$
In particular, among the groups $G_{n,j}$, this condition holds precisely for $j\in\{3,4,5,6,7\}$.
\end{coro}

\begin{proof}
If $\tau\in G$, then $G=H\times\langle\tau\rangle$. Otherwise, projection to $\aut(Q,C_4)$ identifies $G$ with a finite subgroup, and the exponent of $\tau$ defines a homomorphism from this subgroup to $\dz_2$ with kernel $H$.

Finite subgroups of $\pgl_2(\dc)$ are cyclic, dihedral, tetrahedral, octahedral, or icosahedral. A cyclic group has a non-trivial homomorphism to $\dz_2$ precisely when its order is even. For a dihedral group, such a homomorphism can be non-trivial on the reflections, or, when the rotation subgroup has even order, on the rotations. In the latter situation, the two choices for its value on a reflection are conjugate in the normaliser of the rotation subgroup. These possibilities give the seven families in the table.

In case~\ref{joe1}, the only repetitions among these families are
$$
G_{1,2}\sim G_{2,1},\qquad G_{1,5}\sim G_{2,3},\qquad
G_{1,6}\sim G_{1,4},\qquad G_{1,7}\sim G_{2,6},
$$
where $\sim$ denotes conjugacy in $\aut(X)$. The first three follow because $B$ and $R_2$ are conjugate in $\pgl_2(\dc)$; the last follows because the normaliser of $\langle R_2,B\rangle\cong\dz_2^2$ permutes its three subgroups of order two transitively. Apart from these low-order coincidences, the projected subgroup, the presence of $\tau$, and $H$ distinguish the conjugacy classes. Indeed, a homomorphism to $\dz_2$ is uniquely determined by its kernel, so two subgroups not containing $\tau$ are conjugate exactly when their projections and their subgroups $H$ are simultaneously conjugate in $\aut(Q,C_4)$. For a dihedral projection of order greater than four, the cyclic rotation subgroup is characteristic, distinguishing it from kernels containing reflections. In case~\ref{joe2}, conjugation preserves the torus $\dg_m$ and cannot interchange a non-trivial rotation with a reflection, so all seven families are irredundant.

The groups $\mathfrak A_4$ and $\mathfrak A_5$ have no non-trivial homomorphism to $\dz_2$, whereas the only non-trivial homomorphism $\mathfrak S_4\to\dz_2$ is the sign character. The matrices $R_2,B,C$ belong to its kernel, and $D$ does not. This gives $\mathfrak S^{\mathrm{tw}}$.

The remaining cases follow directly from Theorem~\ref{theo:autQC}; in case~\ref{joe3}, every non-trivial finite subgroup of $\dg_a\rtimes\dz_2$ is conjugate to $\langle R_2\rangle\cong\dz_2$. In case~\ref{joe4}, the ambient group is abelian, and its subgroups correspond bijectively to the subspaces $U\subset\mathbb Z_2^3$, so this list is also irredundant.

Finally, $\aut(Q,C_4)$ acts trivially on $\pic(X)$, while $\tau$ exchanges the two extremal rays of the nef cone.
\end{proof}

\subsection{Fixed points on $X$}

In case~\ref{joe1} of Theorem~\ref{theo:autQC}, two points on $Q$ occur both in the superrigidity argument and in the construction of a birational model when $H\cong\mathfrak A_4$ or $H\cong\mathfrak S_4$.

\begin{lemm}\label{lemm:a4-fixed-pair}
Suppose that $Q=\{f_5=0\}$, and let $\mathfrak A=\langle R_2,B,C\rangle\cong\mathfrak A_4$. Its fixed locus on $Q$ consists of the two points
$$
p_\pm=[1:0:\pm\alpha:0:1],\qquad \alpha=\frac{i}{\sqrt 3}.
$$
Both lie outside $C_4$, and we identify them with their lifts to $X$. The normaliser of $\mathfrak A$ in $\aut(Q,C_4)$ is $\mathfrak S=\langle\mathfrak A,D\rangle\cong\mathfrak S_4$. Each element of $\mathfrak S\setminus\mathfrak A$ exchanges $p_+$ and $p_-$, as does $\tau$.
\end{lemm}

\begin{proof}
Put $K=\langle R_2,B\rangle$, the normal Klein four subgroup of $\mathfrak A$. In the coordinates of the Veronese embedding, $R_2$ changes the signs of $x_1,x_3$ and $B$ reverses the coordinates. Their common eigenspaces meet $Q$ only on the line $\{x_1=x_3=0,\ x_0=x_4\}$. On this line the equation of $Q$ is $3x_2^2+x_0^2=0$, so $Q^K=\{p_+,p_-\}$. Neither point belongs to $C_4$. The element $C$ has order three and preserves this pair, so it fixes both points. Hence $Q^{\mathfrak A}=Q^K$.

Since $K$ is characteristic in $\mathfrak A$, its normaliser contains that of $\mathfrak A$. Conjugation on the three involutions of $K$ gives a homomorphism $N_{\pgl_2(\dc)}(K)\to\mathfrak S_3$ whose kernel is $K$ itself. Indeed, an element commuting with $R_2$ preserves or exchanges its two fixed points on $\pl$, and commuting also with $B$ leaves just the four elements of $K$. Thus the normaliser of $\mathfrak A$ has order at most $24$. It contains $\langle\mathfrak A,D\rangle\cong\mathfrak S_4$, so equality holds.

For $t=\pm\alpha$, write $p_t=[1:0:t:0:1]$. Acting in the fourth symmetric power, $D$ sends $p_t$ to
$$
[2-6t:0:-2-2t:0:2-6t]=p_{-t}.
$$
This gives the action of the coset $\mathfrak A D$. Finally, the formula for $\sigma$ in \cite[Example~1.3]{malbon2024automorphisms}, with $s=\tfrac12$, gives $\sigma(p_t)=[-t:0:-\tfrac13:0:-t]=p_{-t}$. Thus $\tau$ exchanges the two lifts on $X$.
\end{proof}

\begin{coro}\label{coro:fixed-points-autX}
Let $G\subset\aut(X)$ be a finite subgroup, written up to conjugacy as in Corollary~\ref{coro:finite-subgroups-autX}.

\begin{enumerate}
\item For the seven families $G_{n,j}$, the following table gives the values of $n$ for which $G_{n,j}$ fixes a point on $X$.

$$
\begin{array}{c|c}
G & \text{Condition for }X^G\neq\emptyset \\ \hline
G_{n,1} & \text{every }n\geq 1\\[2pt]
G_{n,2} & n\in\{1,2\}\\[2pt]
G_{n,3},\ G_{n,4} & \text{every }n\geq 1\\[2pt]
G_{n,5} & n=1\\[2pt]
G_{n,6} & n\in\{1,2,3\}\\[2pt]
G_{n,7} & n\in\{1,2\}.
\end{array}
$$
\item Among the additional groups in case~\ref{joe1}, precisely

$$
\langle R_2,B,C\rangle\cong\mathfrak A_4 \qquad\text{and}\qquad \mathfrak S^{\mathrm{tw}} =\langle R_2,B,C,D\tau\rangle\cong\mathfrak S_4
$$
fix a point on $X$.

\item In cases~\ref{joe3}, \ref{joe5}, and~\ref{joe6}, every finite subgroup fixes a point on $X$.

\item In case~\ref{joe4}, a finite subgroup fixes a point on $X$ if and only if it is a proper subgroup of $\aut(X)\cong\dz_2^3$.
\end{enumerate}
\end{coro}
\begin{proof}
If $\pic^G(X)\not\cong\dz$, then $G=H$. For such groups, the dihedral cases are treated below, and the other cases follow from the Magma script provided in \cite{pinardinwebsite}: whenever $Q^H$ is nonempty, it contains a point outside $C_4$, which lifts to an $H$-fixed point on $X$.

We first show that no subgroup $K\cong\dz_2^3$ fixes a point. Its projection to $\aut(Q,C_4)\subset\pgl_2(\dc)$ cannot be injective, so $\tau\in K$ and $K=L\times\langle\tau\rangle$, where $L=K\cap\aut(Q,C_4)\cong\dz_2^2$. Up to conjugacy, $L=\langle R_2,B\rangle$ in cases~\ref{joe1}, \ref{joe2}, and~\ref{joe4}. Write
$$Q=\{\mu(f_0+f_4)+\lambda f_2+f_5=0\},$$
with $\mu=0$ in cases~\ref{joe1} and~\ref{joe2}, and $\lambda=0$ in case~\ref{joe1}. The common eigenspaces of $R_2$ and $B$ give
$$Q^L=\{p_t:(\lambda+3)t^2-2\mu t+1-\lambda=0\},\qquad p_t=[1:0:t:0:1].$$
These are two distinct points outside $C_4$; the other common eigenspaces miss $Q$ since $\lambda\ne1$ and $\mu\ne\pm2$. The involution formulas in \cite[Example~1.3 and the proof of Theorem~3.1]{malbon2024automorphisms} give
$$\tau:t\longmapsto\frac{\mu-2t}{2-\mu t}=\frac{2\mu}{\lambda+3}-t,$$
which exchanges the two roots. Thus $X^K=\emptyset$. This excludes every group containing $\dz_2^3$, namely $G_{n,5}$ with $n$ even, the three polyhedral groups with $\tau$ adjoined, and the full group in case~\ref{joe4}.

Every cyclic subgroup fixes a point by the holomorphic Lefschetz formula: since $H^i(X,\mathcal O_X)=0$ for $i>0$, its generator has Lefschetz number $1$. Every subgroup isomorphic to $\dz_2^2$ also fixes a point. Indeed, choose an involution $h\in H$. Its fixed locus on $Q$ consists of a smooth conic and two points, and its normal eigenvalues at the fixed points of $C_4$ are $1,-1$. Hence $X^h$ has a unique curve, isomorphic to $\mathbb P^1$, on which the other generator has a fixed point. This settles cases~\ref{joe3}, \ref{joe5}, \ref{joe6}, and all proper subgroups in case~\ref{joe4}.

Let $p_0,\ldots,p_4$ be the coordinate points, and denote by $\ell_{p,q}$ the line passing through two points $p$ and $q$. For $H=\langle R_n,B\rangle\cong D_n$ with $n\ge3$, the eigenspaces of $R_n$ give
$$
Q^{\langle R_n\rangle}=
\begin{cases}
\ell_{p_0,p_3}\cup\ell_{p_1,p_4},&n=3,\\
\{p_0,p_1,p_3,p_4\},&n\ge4.
\end{cases}
$$
Here $\lambda\ne1,-3$; for $n=4$, the additional eigenspace $\langle p_0,p_4\rangle$ meets $Q$ only at $p_0,p_4$. Since $B$ exchanges the two non-intersecting lines $\ell_{p_0,p_3}$ and $\ell_{p_1,p_4}$, or the four points in pairs, $Q^H=\emptyset$. This excludes $G_{n,2},G_{n,5}$ and $G_{n,7}$ for $n\ge3$.

By Corollary~\ref{coro:finite-subgroups-autX}, it remains to consider $G_{n,3},G_{n,6},G_{2,7}$ and $\mathfrak S^{\mathrm{tw}}$. For $G_{2,7}$, the two $H$-fixed points are the points $p_t$ above with $\mu=0$; both $R_4$ and $\tau$ exchange them, so $R_4\tau$ fixes their lifts. For $\mathfrak S^{\mathrm{tw}}=\langle\mathfrak A,D\tau\rangle$, Lemma~\ref{lemm:a4-fixed-pair} shows that $\mathfrak A$ and $D\tau$ fix both $p_+$ and $p_-$.

Use $(\pi,\pi')$ to view $X$ inside $Q\times Q$, where $\tau$ exchanges the factors. For the remaining families, write $s^2=(1-\lambda)/4$, with $s=\tfrac12$ in case~\ref{joe1}, choosing the sign so that
$$
Q=\{s^2x_0x_4-x_1x_3+(1-s^2)x_2^2=0\},\qquad
\sigma=[f_4:sf_3:s^2f_2:sf_1:f_0],
$$
as in \cite[Example~1.3]{malbon2024automorphisms}. This formula shows that $(p_0,p_0)\in X$ is fixed by the torus and by $\tau$, proving the assertion for $G_{n,3}$. For $G_{n,6}$, only $n\ge3$ remains. For $n=3$, the point $[1:0:0:r:0]\in Q\setminus C_4$ is fixed by $R_3$, and $B\sigma$ sends $r$ to $-s/r$; taking $r^2=-s$ gives a $G_{3,6}$-fixed point. For $n\ge4$, the normal weights at $p_0$ are $2,3$, meaning that the differential of $R_n$ on $T_{p_0}Q/T_{p_0}C_4$ has eigenvalues $\zeta_n^2,\zeta_n^3$. Thus there are two $R_n$-fixed points above $p_0$, and the formula for $\sigma$ gives their second coordinates $p_0,p_1$; similarly at $p_4$. Consequently,
$$X^{\langle R_n\rangle}=\{(p_0,p_0),(p_0,p_1),(p_1,p_0),(p_3,p_4),(p_4,p_3),(p_4,p_4)\}.$$
Since $B\tau(p_i,p_j)=(p_{4-j},p_{4-i})$, such a point is fixed by $B\tau$ if and only if $i+j=4$, which none of the listed points satisfies.
\end{proof}

\section{Equivariant rigidity}

The goal of this section is to prove Theorem \ref{theo:rigid}. We will start by recalling a classical result for obstructing equivariant birational maps to $G$-Mori fibre spaces, and apply it together with the results of \cite{pinardin2025A5} to prove the equivariant birational superrigidity of $X$ for $G\cong\ac\times\dz_2$. In all the other cases, we will find an equivariant birational map to a $G$-Mori fibre space over a base of positive dimension.

\begin{prop}\label{prop:notgsol}
Let $G\subset\aut(X)$ be a finite subgroup such that $\pic^G(X)\cong\mathbb Z$ and $G\not\cong\mathfrak A_5\times\mathbb Z_2$. Then $X$ is not $G$-solid. In particular, it is not $G$-rigid.
\end{prop}

\subsection{Noether-Fano inequalities}
We use the discrepancy conventions of \cite[Section~2]{pinardin2025A5}. Let $X$ be a projective variety with at most klt singularities and $G$ a finite subgroup of $\aut(X)$. Let $\mls_X$ be a non-empty mobile $G$-invariant linear system on $X$ consisting of $\q$-Cartier divisors, let $\pi\colon \widetilde X\to X$ be a log-resolution of $(X,\mls_X)$, and
$$
\mls^{\widetilde{X}}:=\pi^*(K_X+\mls_X)-K_{\widetilde X}.
$$
For any prime divisor $D$ on $\widetilde X$, the {\em log discrepancy} of the log pair $(X,\mls_X)$ at $D$ is the rational number
$$
a(X,\mls_X;D)=1-\mathrm{mult}_D(\mls^{\widetilde{X}}),
$$
where the multiplicity is that of a general member of $\mls^{\widetilde{X}}$.

Let $p\in X$ be a point. We say that $(X,\mls_X)$ is {\em canonical (resp. log-canonical, klt)} at $p$ if, for every proper birational morphism $\pi\colon \widetilde X\to X$ with $\widetilde X$ smooth and every prime divisor $E$ on $\widetilde X$ such that $p\in\pi(E)$, we have $a(X,\mls_X;E)\geq 1$ (resp. $\geq 0$, $>0$). The {\em non-canonical (resp. non-log-canonical, non-klt) locus} of $(X,\mls_X)$ is the union of points where $(X,\mls_X)$ is not canonical (resp. not log-canonical, not klt).

An irreducible subvariety $Z\subset X$ is said to be {\em a non-canonical (resp. non-log-canonical, non-klt) centre} of $(X,\mls_X)$ if there exists a proper birational morphism $\pi\colon \widetilde X\to X$ with $\widetilde X$ smooth and a prime divisor $E$ on $\widetilde X$ such that $\pi(E)=Z$ and $a(X,\mls_X;E)<1$ (resp. $<0$, $\leq 0$). To obstruct the existence of $G$-birational maps to $G$-Mori fibre spaces, we will need the classical Noether--Fano inequality, which we recall here in its $G$-equivariant setting as it is in \cite[Theorem 3.2.6]{CS}.

\begin{theo}\label{theo:nfi}
Let $X$ be a Fano variety with terminal singularities, $G$ a finite subgroup of $\aut(X)$ such that $\mathrm{rk}(\cl^G(X))=1$. Assume that there exists a $G$-equivariant birational map $\chi: X\dashrightarrow V$, where $V$ is a variety with a generically free $G$-action such that one of the following  holds:
\begin{enumerate}
\item either $V$ is also a Fano variety with terminal singularities such that\\$\mathrm{rk}(\cl^G(V))=1$;
\item  or there exists a $G$-equivariant morphism $V\to Z$ with connected fibres such that its general fibre is a Fano variety, and $Z$ is a normal projective variety with $\dim(V)>\dim(Z)>0$.
\end{enumerate}
In the former case, let $\mls_X$ be the strict transform on $X$ of the linear system $|-nK_V|$ for $n\gg 0$. In the latter case, let $\mls_X$ be the strict transform on $X$ of the linear system $|H_V|$, where $H_V$ is the pullback on $V$ of a very ample divisor on $Z$ whose class in $\mathrm{Pic}(Z)$ is $G$-invariant. Let $\lambda\in\q$ be such that $\lambda\mls_X\sim_\q-K_X$. Then, if $\chi$ is not biregular, the log pair $(X,\lambda\mls_X)$ has non-canonical singularities.
\end{theo}

\subsection{Superrigidity for $G\cong\mathfrak A_5\times\dz_2$}

The aim of this subsection is to prove the following result.

\begin{prop}\label{prop:A5xZ2rigid}
Let $G\subset\aut(X)$ be a subgroup isomorphic to $\mathfrak A_5\times\dz_2$. Then $\pic^G(X)\cong\dz$ and $X$ is $G$-superrigid.
\end{prop}

We first record the following group-theoretic observation.

\begin{lemm}
Let $G\subset\aut(X)$ be a finite subgroup, and set $H=G\cap\aut(Q,C_4)$. Then $G\cong\mathfrak A_5\times\dz_2$ if and only if $H\cong\mathfrak A_5$ and $\pic^G(X)\cong\dz$. In this case, $G=H\times\langle\tau\rangle$.
\end{lemm}

\begin{proof}
Suppose first that $G\cong\mathfrak A_5\times\dz_2$. Then $[G,G]\cong\mathfrak A_5$ is contained in $\aut(Q,C_4)$. Its centraliser in $\pgl_2(\dc)$ is trivial, so the unique non-trivial central element of $G$ must be $\tau$. Thus $H=[G,G]\cong\mathfrak A_5$ and $G=H\times\langle\tau\rangle$. Since $\tau$ exchanges the two extremal rays of the nef cone of $X$, we obtain $\pic^G(X)\cong\dz$.

Conversely, suppose that $H\cong\mathfrak A_5$ and $\pic^G(X)\cong\dz$. We have $G=H\times\langle\tau\rangle$. Indeed, using the decomposition $\aut(X)=\aut(Q,C_4)\times\langle\tau\rangle$, this follows from the fact that $H$ is self-normalising in $\aut(Q,C_4)\cong\pgl_2(\dc)$, while the condition $\pic^G(X)\cong\dz$ rules out $G=H$.
\end{proof}

For the rest of this subsection, assume that $G\cong\mathfrak A_5\times\dz_2$, and put $H=G\cap\aut(Q,C_4)\cong\mathfrak A_5$. Up to conjugation in $\aut(Q)$, the $H$-action on $Q$ is the one described in \cite[Section~6]{pinardin2025A5}; we transport its notation for invariant curves and orbits to the coordinates fixed above. Let $\mls$ be a non-empty mobile $G$-invariant linear system on $X$, and let $\lambda\in\mathbb Q_{>0}$ satisfy $\lambda\mls\sim_{\mathbb Q}-K_X$. Set $\mls_Q=\pi_*(\mls)$. Since $K_X=\pi^*K_Q+E$, we obtain
\begin{align}
\lambda\mls_Q&\sim_{\mathbb Q}-K_Q,\notag\\ \mult_{C_4}(\lambda\mls_Q)&=1,\notag\\ K_X+\lambda\mls&=\pi^*(K_Q+\lambda\mls_Q).\label{eq:discrep}
\end{align}

Let $Z_0$ be a non-canonical centre of $(X,\lambda\mls)$, and let $Z=G\cdot Z_0$ be the support of its $G$-orbit. By \eqref{eq:discrep}, every irreducible component of $\pi(Z)$ is a non-canonical centre of $(Q,\lambda\mls_Q)$. Thus, by \cite[Proposition~6.1]{pinardin2025A5}, each $H$-irreducible component $\gamma$ of $\pi(Z)$ satisfies
$$
\gamma\in\{\Sigma_5,\Sigma_5',C_4,C_4',C_8,C_8'\},
$$
or is a smooth irreducible curve of degree $10$ and genus $6$. Here $\Sigma_5$ and $\Sigma_5'$ are the two $H$-orbits of length $5$ on $Q$, the curves $C_4$ and $C_4'$ are the two $H$-invariant rational normal quartics, and $C_8$ and $C_8'$ are the two $H$-invariant smooth rational curves of degree $8$.

\begin{rema}\label{rema:bound-deg}
Let $C\subset X$ be a curve that is a non-canonical centre of $(X,\lambda\mls)$. For general members $M_1,M_2\in\lambda\mls$, we have

$$
M_1\cdot M_2=mC+\Delta,
$$
for some $m>1$ and some effective $1$-cycle $\Delta$ not containing $C$ in its support. Intersecting with $-K_X$, we obtain
\begin{align*}
-K_X\cdot C &<-K_X\cdot(mC)\\ &=-K_X\cdot(M_1\cdot M_2-\Delta)\\ &\le -K_X\cdot M_1\cdot M_2\\ &=(-K_X)^3=28.
\end{align*}
The same argument applies to a $G$-irreducible curve whose irreducible components are non-canonical centres, since the coefficients of its components in $M_1\cdot M_2$ are equal.
\end{rema}

\begin{lemm}\label{lemm:swap-five-orbits}
The $H$-orbits $\Sigma_5$ and $\Sigma_5'$ on $Q$ are exchanged by the birational involution $\sigma$.
\end{lemm}

\begin{proof}
Fix a subgroup $\mathfrak A\subset H$ isomorphic to $\mathfrak A_4$. Each orbit of length five is isomorphic to $H/\mathfrak A$ and contains exactly one $\mathfrak A$-fixed point, since $\mathfrak A$ is self-normalising in $H\cong\mathfrak A_5$. After conjugating $\mathfrak A$ in $\aut(Q,C_4)$, these two points are $p_+$ and $p_-$ of Lemma~\ref{lemm:a4-fixed-pair}. The involution $\sigma$ exchanges them and is $H$-equivariant, so it exchanges their $H$-orbits.
\end{proof}

\begin{lemm}\label{lemm:quadrics-five-orbits}
The base locus of the linear system $\mathcal Q_{\Sigma_5\cup\Sigma_5'}$ of quadric sections of $Q$ containing $\Sigma_5\cup\Sigma_5'$ is exactly $\Sigma_5\cup\Sigma_5'$.
\end{lemm}

\begin{proof}
We use the points $p_+,p_-$ of Lemma~\ref{lemm:a4-fixed-pair} and the Magma script provided in \cite{pinardinwebsite}. Put
$$
\varphi=\frac{1+\sqrt5}{2}, \qquad
g=\frac12\begin{pmatrix}-\varphi+i&\varphi-1\\\varphi-1&\varphi+i\end{pmatrix}.
$$
Direct multiplication in $\pgl_2(\dc)$ gives $R_2^2=g^3=(R_2g)^5=\mathrm{id}$, and the script checks that $C\in\langle R_2,g\rangle$. By the presentation $\mathfrak A_5\cong\langle a,g\mid a^2=g^3=(ag)^5=1\rangle$ in \cite[Presentations]{atlasA5}, we obtain $\langle R_2,g\rangle\cong\mathfrak A_5$, containing $\langle R_2,C\rangle=\mathfrak A$. After conjugating $H$, we may thus take $H=\langle R_2,g\rangle$ and, up to interchanging the two orbits, write
$$
\Sigma_5=\{(R_2g)^j(p_+):0\le j\le4\}, \qquad
\Sigma_5'=\{(R_2g)^j(p_-):0\le j\le4\},
$$
where $R_2g$ acts through the fourth symmetric power.

The script computes a basis of the linear system and verifies that its base scheme is zero-dimensional, reduced, and of degree $10$. Since it contains the ten distinct points of $\Sigma_5\cup\Sigma_5'$, it is precisely their union.
\end{proof}

\begin{lemm}\label{lemm:exclude-five-orbits}
No point of $\pi^{-1}(\Sigma_5\cup\Sigma_5')$ is a non-canonical centre of $(X,\lambda\mls)$.
\end{lemm}

\begin{proof}
Set $\Sigma=\Sigma_5\cup\Sigma_5'$. Since point stabilisers for the $H$-action on $C_4\cong\pl$ are cyclic of order at most $5$, every orbit has length at least $12$, so we have $\Sigma\cap C_4=\emptyset$. Thus the morphism $\pi$ is an isomorphism in a neighbourhood of each point of $\pi^{-1}(\Sigma)$. By Lemma~\ref{lemm:swap-five-orbits}, the set $\pi^{-1}(\Sigma)$ is a single $G$-orbit. Suppose that its points are non-canonical centres of $(X,\lambda\mls)$. Then, by \eqref{eq:discrep}, the points of $\Sigma$ are non-canonical centres of $(Q,\lambda\mls_Q)$.

The stabiliser in $H$ of any point $p\in\Sigma$ is isomorphic to $\aq$ and acts faithfully on the tangent space $T_pQ$. The irreducible complex representations of $\aq$ have dimensions $1,1,1,3$, and every one-dimensional representation is trivial on its normal Klein four subgroup; see, for example, \cite{dokchitserGroupNames}. Thus this tangent representation is irreducible. By \cite[Lemma~8.3]{pinardin2025A5}, we obtain

$$
\mult_p(\lambda\mls_Q)>2,\forall p\in\Sigma.
$$
Let $M_1,M_2$ be general members of $\lambda\mls_Q$, and let $S$ be a general member of $\mathcal Q_{\Sigma}$. By Lemma~\ref{lemm:quadrics-five-orbits}, this linear system has finite base locus, so $S$ contains no irreducible component of $M_1\cap M_2$. We obtain
\begin{align*}
36&=S\cdot M_1\cdot M_2\\ &\ge\sum_{p\in\Sigma}\mult_p(M_1)\mult_p(M_2)\\ &=\sum_{p\in\Sigma}\mult_p(\lambda\mls_Q)^2>40,
\end{align*}
a contradiction.
\end{proof}

\begin{lemm}\label{lemm:centers-in-E}
Let $Z$ be a $G$-invariant subvariety of $X$. If $Z$ contains a non-canonical centre of $(X,\lambda\mls)$, then $Z$ is not contained in $E\cup\tau(E)$.
\end{lemm}

\begin{proof}
Suppose that $Z\subset E\cup\tau(E)$, and let $W\subset Z$ be a non-canonical centre of $(X,\lambda\mls)$. Since $\tau\in G$, replacing $W$ by $\tau(W)$ if necessary, we may assume that $W\subset E$. By \eqref{eq:discrep}, its image $\pi(W)\subset C_4$ is a non-canonical centre of $(Q,\lambda\mls_Q)$. We have $\mult_{C_4}(\lambda\mls_Q)=1$, so $C_4$ itself is not a non-canonical centre. Therefore $\pi(W)$ is a point. Its $H$-orbit has length at least $12$, contradicting \cite[Proposition~6.1]{pinardin2025A5}.
\end{proof}

Let $R$ and $R'$ be the two non-normal cubic sections of $Q$ defined in \cite[(6.4)]{pinardin2025A5}. Their singular loci are $C_4$ and $C_4'$, respectively.

\begin{lemm}\label{lemm:exclude-degree-ten}
Let $C\subset X$ be an irreducible $H$-invariant curve such that $\pi(C)$ has degree $10$. Then $C$ is not a non-canonical centre of $(X,\lambda\mls)$.
\end{lemm}

\begin{proof}
Suppose that $C$ is a non-canonical centre, and set $\gamma=\pi(C)$. By \eqref{eq:discrep}, the curve $\gamma$ is a non-canonical centre of $(Q,\lambda\mls_Q)$. By \cite[Lemma~6.17]{pinardin2025A5}, it is not contained in $R\cup R'$.

The proof of \cite[Lemma~6.9]{pinardin2025A5} shows that $\gamma\cdot R=30$, with intersection multiplicity one at each of the thirty intersection points. Since $R$ is singular along $C_4$, any point of $\gamma\cap C_4$ would have intersection multiplicity at least two in $\gamma\cap R$. Hence $\gamma\cap C_4=\emptyset$, and therefore $C\cap E=\emptyset$. We deduce that

$$
-K_X\cdot C=\bigl(\pi^*(-K_Q)-E\bigr)\cdot C=3\deg(\gamma)=30,
$$
contradicting Remark~\ref{rema:bound-deg}.
\end{proof}

\begin{lemm}\label{lemm:exclude-C8-prime}
Let $C\subset X$ be an irreducible curve such that $\pi(C)=C_8'$. Then $C$ is not a non-canonical centre of $(X,\lambda\mls)$.
\end{lemm}

\begin{proof}
Suppose that $C$ is a non-canonical centre. Since $\tau\in G$, the curve $\tau(C)$ is also a non-canonical centre of $(X,\lambda\mls)$. By \cite[Remark~6.16]{pinardin2025A5}, its image $\pi(\tau(C))=\sigma(C_8')$ is an irreducible curve of degree $16$. By \eqref{eq:discrep}, this image is a non-canonical centre of $(Q,\lambda\mls_Q)$, contradicting \cite[Proposition~6.1]{pinardin2025A5}.
\end{proof}

The only centres that are not excluded from the possible non-canonical singularities of $X$ are the strict transforms of $C_4'$ and $C_8$. Recall from \cite[Lemma~6.18]{pinardin2025A5} that $\sigma(C_4')=C_8$, and let $Z_{4,8}$ be the $G$-irreducible curve $\widetilde C_4'\cup\widetilde C_8.$ To prove Proposition~\ref{prop:A5xZ2rigid}, it remains to prove the following.

\begin{lemm}\label{lemm:exclude-Z48}
The irreducible components of $Z_{4,8}$ cannot be non-canonical centres of $(X,\lambda\mls)$.
\end{lemm}

\begin{proof}
Suppose otherwise. Since $\tau$ exchanges $\widetilde C_4'$ and $\widetilde C_8$, and $\pi$ is an isomorphism at their generic points, we have

$$
m=\mult_{C_4'}(\lambda\mls_Q)=\mult_{C_8}(\lambda\mls_Q)>1.
$$
Let $\sigma'\colon Q\dashrightarrow Q$ be the $H$-equivariant birational involution defined by the quadrics through $C_4'$, as in \cite[Lemma~6.15 and Remark~6.16]{pinardin2025A5}. Let $\mls_Q'=(\sigma')_*(\mls_Q).$ We first compute the class of $\lambda\mls_Q'$. Let $\rho\colon Y\to Q$ be the blow-up of $C_4'$, let $F$ be its exceptional divisor, and put $L=\rho^*{H_Q}$. The involution $\sigma'$ lifts to an automorphism $\tau'$ of $Y$ satisfying $\tau'^*L\sim 2L-F.$ Since $-K_Y\sim 3L-F$ is preserved by $\tau'$, we also have $\tau'^*F\sim 3L-2F.$ The strict transform $\mls_Y$ of $\mls_Q$ satisfies $\lambda\mls_Y\sim_{\mathbb Q}3L-mF.$ Hence
\begin{align*}
\lambda\tau'_*(\mls_Y)&\sim_{\mathbb Q}3(2L-F)-m(3L-2F)\\ &=(6-3m)L-(3-2m)F.
\end{align*}
Pushing forward by $\rho$, we obtain $\lambda\mls_Q'\sim_{\mathbb Q}3(2-m){H_Q},$ where $H_Q=\mathcal O_Q(1).$ Since $\mls_Q'$ is a non-empty mobile linear system of positive degree, we have $m<2$. Thus, setting $\lambda'=\frac{\lambda}{2-m},$ we obtain $\lambda'>\lambda, \lambda'\mls_Q'\sim_{\mathbb Q}-K_Q.$ By \cite[Remark~6.16]{pinardin2025A5}, the curve $\Gamma=\sigma'(C_8)$ is irreducible of degree $16$. The involution $\sigma'$ is an isomorphism at the generic point of $C_8$: its indeterminacy locus is $C_4'$, and its exceptional divisor $R'$ is contracted onto $C_4'$, while $C_8$ is not contained in $R'$ by \cite[Lemma~6.17]{pinardin2025A5}. Consequently,

$$
\mult_\Gamma(\lambda'\mls_Q')=\frac{\mult_{C_8}(\lambda\mls_Q)}{2-m}=\frac{m}{2-m}>1.
$$
Therefore $\Gamma$ is an $H$-invariant non-canonical centre of $(Q,\lambda'\mls_Q')$. Since $\mls_Q'$ is mobile and $H$-invariant, this contradicts \cite[Proposition~6.1]{pinardin2025A5}.
\end{proof}

\begin{proof}[Proof of Proposition \ref{prop:A5xZ2rigid}]
Let $\mls_X$ be a non-empty mobile $G$-invariant linear system on $X$, and $\lambda\in\q$ such that $\lambda\mls_X\sim_\q-K_X$. Combining the preceding lemmas, there cannot be a non-canonical centre of $(X,\lambda\mls_X)$, and we conclude using Theorem \ref{theo:nfi}.
\end{proof}

\subsection{The cases $H\cong\mathfrak A_4$ and $H\cong\mathfrak S_4$}

We first prove Proposition~\ref{prop:notgsol} for these two possibilities for $H=G\cap\aut(Q,C_4)$. The two points of Lemma~\ref{lemm:a4-fixed-pair} determine a common birational model.
\begin{lemm}\label{lemm:a4s4-actions}
Let $G\subset\aut(X)$ be a finite subgroup such that $\pic^G(X)\cong\dz$, and set $H=G\cap\aut(Q,C_4)$. Suppose that $H\cong\aq$ or $H\cong\sq$. Let $\mathfrak A$ be the unique subgroup of $H$ isomorphic to $\aq$, and put
$$\mathfrak S=N_{\aut(Q,C_4)}(\mathfrak A)\cong\sq.$$
Let $\epsilon\colon\mathfrak S\to\dz_2$ be the homomorphism with kernel $\mathfrak A$.

Identify the two $\mathfrak A$-fixed points with $p_+,p_-$ of Lemma~\ref{lemm:a4-fixed-pair} after conjugation. Exactly one of the following occurs.
\begin{enumerate}
\item $H=\mathfrak A$ and $G=H\times\langle\tau\rangle$. The group $H$ fixes both points and $\tau$ exchanges them.\label{a4s4-1}
\item $H=\mathfrak A$ and $G=\mathfrak S^{\mathrm{tw}}=\{s\tau^{\epsilon(s)}\mid s\in\mathfrak S\}$. The group $G$ fixes both points.\label{a4s4-2}
\item $H=\mathfrak S$ and $G=\mathfrak S\times\langle\tau\rangle$. Both $\tau$ and every element of $\mathfrak S\setminus\mathfrak A$ exchange the two points.\label{a4s4-3}
\end{enumerate}
\end{lemm}

\begin{proof}
Corollary~\ref{coro:finite-subgroups-autX} gives the three possibilities for $G$. By Lemma~\ref{lemm:a4-fixed-pair}, an element $s\tau^e$ with $s\in\mathfrak S$ exchanges the points precisely when $\epsilon(s)+e=1$ in $\dz_2$. This gives the asserted action in each case.
\end{proof}

\begin{prop}\label{prop:a4s4-rigid}
Let $G$ be as in Lemma~\ref{lemm:a4s4-actions}. There is a $G$-equivariant birational map
$$
X\dashrightarrow W\subset\mathbb P^2\times\mathbb P^2,
$$
where $W$ is a smooth divisor of bidegree $(1,1)$. The action of $G$ on the two factors is its permutation action on $\{p_+,p_-\}$.
\end{prop}

\begin{proof}
Let $p_+$ and $p_-$ be the two $\mathfrak A$-fixed points from Lemma~\ref{lemm:a4s4-actions}, and let $\rho_\pm\colon Q\dashrightarrow\mathbb P^3$ be the projections from $p_\pm$. Denote by $C_{4,\pm}$ the image of $C_4$ under $\rho_\pm$, and put $q_\pm=\rho_\pm(p_\mp)$. Define

$$
V_\pm =H^0\!\left( \mathbb P^3, \mathcal I_{C_{4,\pm}} \mathcal I_{q_\pm}^{\,2}(3) \right).
$$
Thus $V_\pm$ is the vector space of cubic forms vanishing along $C_{4,\pm}$ and to order at least two at $q_\pm$. Using the Magma script provided in \cite{pinardinwebsite}, we compute $\dim V_+=\dim V_-=3$. These spaces therefore define rational maps

$$
\eta_\pm\colon \mathbb P^3\dashrightarrow \mathbb P(V_\pm^\vee)\cong\mathbb P^2.
$$
Consider the rational map

$$
\psi =(\eta_+\circ\rho_+,\eta_-\circ\rho_-) \colon Q\dashrightarrow\mathbb P^2\times\mathbb P^2.
$$
The same script computes the closure $W$ of its image and shows that $\psi$ is birational onto $W$. After choosing suitable coordinates on the two factors, we obtain

$$
W=\{x_1x_2+y_1y_2+z_1z_2=0\} \subset\mathbb P^2\times\mathbb P^2.
$$
Set $\phi=\psi\circ\pi\colon X\dashrightarrow W$.

To establish equivariance, let $\beta\colon Y\to X$ be the blow-up of $p_+$ and $p_-$, with exceptional divisors $F_+$ and $F_-$. Since $\{p_-,p_+\}$ is $G$-invariant, the $G$-action lifts to $Y$. The two components of $\phi\circ\beta$ are defined by the complete linear systems

$$
\mathcal L_\pm=|-K_Y-F_\pm|.
$$
Indeed, writing $D=(\pi\circ\beta)^*\mathcal O_Q(1)$ and denoting by $E_Y$ the strict transform of the $\pi$-exceptional divisor, we have

$$
-K_Y-F_\pm =3D-E_Y-3F_\pm-2F_\mp.
$$
These systems correspond precisely to cubic sections of $Q$ containing $C_4$, with multiplicity at least three at $p_\pm$ and at least two at $p_\mp$. Projection from $p_\pm$ identifies them with $\mathbb P(V_\pm)$.

It follows that $G$ permutes $\mathcal L_+$ and $\mathcal L_-$ exactly as it permutes $p_+$ and $p_-$. Hence $G$ acts biregularly on $\mathbb P(V_+^\vee)\times\mathbb P(V_-^\vee)$, preserving $W$, and $\phi$ is $G$-equivariant. An element fixing both points preserves both factors, whereas an element interchanging the points interchanges the factors.

Let $L_1,L_2$ denote the restrictions to $W$ of the hyperplane classes from the two factors of $\pp\times\pp$. Then $\pic(W)=\dz L_1\oplus\dz L_2,$ and $-K_W=2L_1+2L_2.$ In particular, $W$ is Fano, and both projections $W\to\mathbb P^2$ are $\mathbb P^1$-bundles.

\end{proof}

\begin{prop}\label{prop:a4s4-not-solid}
Let $G\subset\aut(X)$ be a finite subgroup and put $H=G\cap\aut(Q,C_4)$. If $H\cong\mathfrak A_4$ or $H\cong\mathfrak S_4$, then $X$ is not $G$-solid.
\end{prop}

\begin{proof}
Suppose first that $\pic^G(X)\cong\dz$, and consider the model $W$ of Proposition~\ref{prop:a4s4-rigid}. For $G=\mathfrak S^{\mathrm{tw}}$, both projections $W\to\mathbb P^2$ are $G$-equivariant $\mathbb P^1$-bundles of relative Picard rank one. Either projection gives the required $G$-Mori fibre space.

Now let $G=\mathfrak A\times\langle\tau\rangle$ or $G=\mathfrak S\times\langle\tau\rangle$. The involution $\tau$ exchanges the factors, so $\pic^G(W)=\dz(L_1+L_2)$, where $L_i$ is the hyperplane class from the $i$-th factor. The model $W$ is thus a $G$-Fano threefold distinct from $X$ (their indices are two and one, respectively). To obtain a fibration over a positive-dimensional base, we construct a $G$-invariant pencil on $W$.

Let $K\cong\dz_2^2$ be the normal Klein four subgroup of $\mathfrak A$. Then $K\triangleleft G$, and $\mathfrak A$ acts faithfully on each factor. Indeed, the action on one factor determines the action on its dual via the incidence relation defining $W$.

Every subgroup isomorphic to $\dz_2^2$ in $\operatorname{PGL}_3(\dc)$ is conjugate to the diagonal subgroup of order four; see \cite[Section~3.2]{pinardin2024linearization}. Indeed, its inverse image in $\operatorname{SL}_3(\dc)$ has a Sylow $2$-subgroup mapping isomorphically onto it, so its action can be simultaneously diagonalised. Thus, choosing suitable coordinates on the first factor and dual coordinates on the second, we may write
$$
W=\{x_0y_0+x_1y_1+x_2y_2=0\}\subset\mathbb P^2_{[x_0:x_1:x_2]}\times\mathbb P^2_{[y_0:y_1:y_2]},
$$
with $K$ generated by $(x,y)\mapsto(D_1x,D_1y)$ and $(x,y)\mapsto(D_2x,D_2y)$, where
$$
D_1=\operatorname{diag}(1,-1,-1),\qquad D_2=\operatorname{diag}(-1,1,-1).
$$
Since $K\triangleleft G$, every element of $G$ preserving the two factors permutes the three common eigenlines of $K$. It therefore acts by a monomial matrix on the first factor and its inverse transpose on the second. Such an element permutes the products $x_0y_0,x_1y_1,x_2y_2$. The same conclusion holds for elements exchanging the factors, by composing with $(x,y)\mapsto(y,x)$. Consequently, the pencil generated by these products is $G$-invariant and defines a $G$-equivariant dominant rational map
$$
f:W\dashrightarrow B=\{t_0+t_1+t_2=0\}\cong\mathbb P^1,\qquad (x,y)\longmapsto[x_0y_0:x_1y_1:x_2y_2].
$$

\noindent For $t=[t_0:t_1:t_2]\in B$ with $t_0t_1t_2\ne0$, the closure of the fibre is the surface
$$
S_t=\{t_1x_0y_0-t_0x_1y_1=t_2x_0y_0-t_0x_2y_2=0\}\subset\mathbb P^2\times\mathbb P^2.
$$
This surface is the closure of the graph of the quadratic transformation
$$
\mathbb P^2\dashrightarrow\mathbb P^2,\qquad [x_0:x_1:x_2]\longmapsto[t_0x_1x_2:t_1x_0x_2:t_2x_0x_1].
$$
Hence $S_t$ is isomorphic to the blow-up of $\mathbb P^2$ at its three coordinate points. In particular, it is a smooth del Pezzo surface of degree $6$.

Resolve the indeterminacy of $f$ equivariantly to obtain a morphism $\widetilde W\to B$. Its general fibre is a smooth rational surface, whose canonical class is not in the closure of the cone of effective divisor classes. Thus $K_{\widetilde W}$ is not pseudoeffective over $B$. The relative $G$-equivariant minimal model program therefore terminates with a $G$-Mori fibre space $Y\to Z$, together with a dominant morphism $Z\to B$. Thus $\dim Z\ge1$. Since $Y$ is $G$-birational to $X$, the threefold $X$ is not $G$-solid.

Finally, if $\pic^G(X)$ has rank two, then $G=H$. Apply the construction above to $H\times\langle\tau\rangle$. Its pencil on $W$ is also $G$-invariant, with the same rational general fibre. Resolving the pencil and running the relative $G$-MMP again gives a $G$-Mori fibre space over a positive-dimensional base.
\end{proof}

\subsection{Non-solidity for the remaining groups}

For the groups $G_{n,j}$ we use the rational quotient by $\dg_m$. The additive case admits an analogous quotient by $\dg_a$. For the three remaining types of quadric, we construct equivariant pencils directly.

\begin{lemm}\label{lemm:gmz2z2}
Let $G\subset\aut(X)$ be isomorphic to $\left(\dg_m\rtimes\dz_2\right)\times\dz_2$. Then $X$ is not $G$-solid.
\end{lemm}
\begin{proof}
Up to a choice of homogeneous coordinates $[x_0:x_1:x_2:x_3:x_4]\in\p^4$ compatible with the Veronese embedding
$$
\mathcal V\colon [u:v]\longmapsto [u^4:u^3v:u^2v^2:uv^3:v^4],
$$
the quadric $Q$ is in case \ref{joe1} or \ref{joe2} of Theorem \ref{theo:autQC}, and the subgroup $\dg_m\subset\aut(Q,C_4)$ is induced by $t\cdot [u:v]=[u:tv],$ hence acts on $\p^4$ by

$$
t\cdot [x_0:x_1:x_2:x_3:x_4]=[x_0:t x_1:t^2 x_2:t^3 x_3:t^4 x_4].
$$
The $\mathbb G_m$-invariant part of the field of functions of $Q$ is $\mathbb{C}(Q)^{\mathbb{G}_m}=\mathbb{C}\left(\frac{x_0 x_4}{x_2^2},\frac{x_0 x_2}{x_1^2}\right).$ Let $\phi$ be the corresponding rational quotient map $Q\dashrightarrow\pp$. Then $\phi\circ\pi$ is a $\dg_m$-invariant dominant rational map $X\dashrightarrow \pp$. Moreover, since $\mathbb G_m$ is normal in $G$, the $\mathbb G_m$-orbits are permuted by the action of $G$, so that it induces a birational action of $G$ on $\mathbb P^2$. Taking a regularisation of this action, we get a $G$-equivariant dominant rational map $\psi\colon X\dashrightarrow B$, where $B$ is a smooth rational surface. Let
$$
\begin{tikzcd} & W & \\ X && B \arrow[from=1-2, to=2-1] \arrow["{\widetilde\psi}", from=1-2, to=2-3] \arrow["\psi"', dashed, from=2-1, to=2-3] \end{tikzcd}
$$
be a $G$-equivariant resolution of indeterminacies, so that $\tilde\psi\colon W\to B$ is a $G$-equivariant fibration. General fibres of $\tilde\psi$ are the closures of $\mathbb G_m$-orbits, so a general fibre of $\tilde\psi$ is a smooth rational curve. Running a $G$-equivariant relative MMP over $B$, we obtain a minimal $G$-conic bundle $Y$ over a rational surface $B'$ such that $Y$ is $G$-birational to $X$.
\end{proof}

\begin{coro}\label{coro:gmsubgroups}
Let $G\subset\aut(X)$ be one of the groups $G_{n,j}$ in cases~\ref{joe1} and~\ref{joe2} of Theorem~\ref{theo:autQC}. Then $X$ is not $G$-solid.
\end{coro}
\begin{proof}
The rational quotient constructed in the proof of Lemma~\ref{lemm:gmz2z2} remains $G$-equivariant after restriction and has rational general fibres. Resolving it $G$-equivariantly and running the relative $G$-MMP yields a $G$-conic bundle over a rational surface.
\end{proof}

\begin{lemm}\label{lemm:gaz2z2}
Let $G\subset\aut(X)$ be isomorphic to $\left(\dg_a\rtimes\dz_2\right)\times\dz_2$. Then $X$ is not $G$-solid.
\end{lemm}
\begin{proof}
We proceed exactly as in the proof of Lemma \ref{lemm:gmz2z2}. The quadric $Q$ is in case \ref{joe1} or \ref{joe3} of Theorem \ref{theo:autQC}, and the action of $\dg_a$ on $C_4$ is induced by $t\cdot [u:v]=[u+tv:v]$. Hence, the $\mathbb G_a$-action on $Q$ is given by

$$
\begin{tikzcd} {\left(t,[x_0:x_1:x_2:x_3:x_4]\right)} \\ {\left[\sum_{i=0}^{4}\binom{4}{i}t^ix_i:\sum_{i=0}^{3}\binom{3}{i}t^ix_{i+1}:\sum_{i=0}^{2}\binom{2}{i}t^ix_{i+2}:\sum_{i=0}^{1}\binom{1}{i}t^ix_{i+3}:x_4\right]} \arrow[maps to, from=1-1, to=2-1] \end{tikzcd}
$$
On the affine chart $\{x_4\neq 0\}$, put $y_i=\frac{x_i}{x_4},$ for $i\in\{0,1,2,3\}.$ Then the $\dg_a$-action is given by
\begin{align*}
y_3&\longmapsto y_3+t,\\ y_2&\longmapsto y_2+2y_3t+t^2,\\ y_1&\longmapsto y_1+3y_2t+3y_3t^2+t^3,\\ y_0&\longmapsto y_0+4y_1t+6y_2t^2+4y_3t^3+t^4.
\end{align*}
We want to compute the $\dg_a$-invariant part of $\dc(Q)$. Take $P_1=y_2-y_3^2,$ and $P_2=y_1-3y_3y_2+2y_3^3$. These two polynomials are $\mathbb G_a$-invariant. The equation of $Q$ expresses $y_0$ as a polynomial in $y_3,y_2,y_1$. Moreover, $y_2=P_1+y_3^2,$ and $y_1=P_2+3y_3P_1+y_3^3,$ so $\mathbb C(Q)=\mathbb C(y_3,y_2,y_1)=\mathbb C(y_3,P_1,P_2).$ Since the $\dg_a$-action fixes $P_1,P_2$ and sends $y_3$ to $y_3+t$, we obtain $\mathbb C(Q)^{\mathbb G_a}=\mathbb C(P_1,P_2),$ which induces a $\dg_a$-invariant dominant rational map $\phi\colon Q\dashrightarrow \p^2$. The rest of the proof is similar to that of Lemma \ref{lemm:gmz2z2}.
\end{proof}

\begin{lemm}\label{lemm:Gis22or222}
Assume that $\aut(X)$ is isomorphic to $\dz_2^2$ or $\dz_2^3$. Then $X$ is not $G$-solid.
\end{lemm}
\begin{proof}
The quadric $Q$ is in case \ref{joe4} or \ref{joe5} of Theorem \ref{theo:autQC}. The conic $\gamma=Q\cap\{x_1=x_3=0\}$ is invariant under $H$, and so is the pencil $\mathcal P$ of all hyperplane sections of $Q$ containing $\gamma$. Let $\mathcal P'=\sigma(\mathcal P)$. By \cite[Remark~1.2, Theorem~3.1]{malbon2024automorphisms}, we have $\tau^*L=2L-E$, where $L=\pi^*\mathcal O_Q(1)$ and $E$ is the exceptional divisor. Since $\tau$ commutes with $g_1$, its pullback identifies the $-1$-eigenspaces of $g_1$ in $H^0(X,L)$ and in $H^0(X,2L-E)$, namely $\tau^*\langle\widetilde x_1,\widetilde x_3\rangle=\langle\widetilde f_1,\widetilde f_3\rangle$. Thus $\mathcal P'$ is the pencil generated by $f_1$ and $f_3$.\\

\noindent For $a\in\dc$, put $\ell_a=\{x_2=ax_0,\ x_3=ax_1,\ x_4=a^2x_0\}\subset\mathbb P^4$. A direct computation, included in the Magma script provided in \cite{pinardinwebsite}, shows that the base locus of $\mathcal P'$ is $\gamma\cup C_4\cup\ell_{a_1}\cup\ell_{a_2}$, where $a_1,a_2$ are the two roots of $\mu a^2-4a+\mu=0$ in case \ref{joe4}, and $\{a_1,a_2\}=\{0,4\}$ in case \ref{joe5}.\\

\noindent Choose a general member $S\in\mathcal P$, which is a smooth quadric surface. No member $S'\in\mathcal P'$ contains $S$ since $S'$ contains the non-degenerate curve $C_4$. Hence, the restriction of $\mathcal P'$ to $S$ defines a pencil. Identify $S$ with $\pl\times\pl$. For $S'\in\mathcal P'$, write $S'|_S=\gamma+\Gamma_{S'}.$ The divisors $S'|_S$ and $\gamma$ have bidegrees $(2,2)$ and $(1,1)$, respectively, so $\Gamma_{S'}$ has bidegree $(1,1)$. Moreover, $S$ meets $C_4$ in four distinct points, exactly two of which, denoted $p,q$, lie outside $\gamma$. Hence $\Gamma_{S'}$ contains $p$ and $q$ for any $S'\in\mathcal P'$. We deduce that the pencil $\{\Gamma_{S'},S'\in\mathcal P'\}$ is the pencil of $(1,1)$-curves passing through $p$ and $q$.\\

\noindent The points $p,q$ do not lie on the same ruling fibre of $S$. Indeed, suppose that they lie on a ruling line $R$. Every residual curve $\Gamma_{S'}$ contains $p,q$ and satisfies $\Gamma_{S'}\cdot R=1$, so it must contain $R$. Consequently, $R$ belongs to the base locus of $\mathcal P'$, and must be one of the lines $\ell_{a_1},\ell_{a_2}$. This is impossible because a general member $S\in\mathcal P$ contains neither of these lines. Thus a general curve $\Gamma_{S'}$ is smooth and isomorphic to $\pl$.\\

\noindent Let $B$ be the union of the base loci of $\mathcal P$ and $\mathcal P'$. For general $S,S'$, the residual conic $\Gamma_{S'}$ meets $B$ in finitely many points. Hence we have $\overline{(S\cap S')\setminus B}=\Gamma_{S'},$ which is a smooth rational curve.\\

\noindent Finally, let $\xi\colon Q\dashrightarrow\pl$ be the rational map defined by $\mathcal P$, and set
\begin{equation}
\phi=(\xi,\xi\circ\sigma)\colon Q\dashrightarrow\pl\times\pl,\qquad \psi=\phi\circ\pi.
\end{equation}
The restriction of $\xi\circ\sigma$ to a general member $S\in\mathcal P$ is nonconstant, so $\phi$, and hence $\psi$, is dominant. Since $\tau$ commutes with $H$, the map $\psi=(\xi\circ\pi,\xi\circ\pi\circ\tau)$ is $G$-equivariant, where $H$ acts diagonally on $\pl\times\pl$ and $\tau$ exchanges the factors.\\

\noindent The closure of a general fibre of $\phi$ is the smooth conic $\Gamma_{S'}$. Resolve $\psi$ equivariantly to obtain a morphism to $\pl\times\pl$. A general fibre of this morphism is a smooth projective curve birational to $\Gamma_{S'}$, hence isomorphic to $\pl$. The relative $G$-MMP therefore gives a $G$-Mori fibre space over a surface.
\end{proof}

\begin{lemm}\label{lemm:Gis2}
If $\aut(X)\cong\dz_2$, then $X$ is not $G$-solid.
\end{lemm}
\begin{proof}
The quadric $Q$ is in case \ref{joe6} of Theorem \ref{theo:autQC}. Here, $Q=\{f_1+f_5=0\}$. We use the same construction as in the preceding proof, this time with $\mathcal P=\{ax_3+bx_4=0\}.$ By the proof of \cite[Theorem~3.1, case~(3)]{malbon2024automorphisms}, we have $\mathcal P'=\sigma(\mathcal P)=\{cf_0+df_1=0\}.$ Let $\ell=\{x_2=x_3=x_4=0\}$. The base loci of $\mathcal P$ and $\mathcal P'$ have supports $\ell$ and $\ell\cup C_4$, respectively. For general $a,b\in\dc$, denote $S_a=Q\cap\{x_3=ax_4\},$ and $S'_b=Q\cap\{f_1=bf_0\}$. The surface $S_a$ is a quadric cone with vertex $v_a=[1+4a:1:0:0:0].$ We have $f_1-bf_0=x_4\left((a+b)x_2-x_1-ba^2x_4\right).$ Moreover, the hyperplane section $\{x_4=0\}$ of $S_a$ is the double line $2\ell$. Consequently, $S'_b|_{S_a}=2\ell+\Gamma_{a,b}$, where $\Gamma_{a,b}=S_a\cap\{x_1=(a+b)x_2-ba^2x_4\}.$ The plane defining $\Gamma_{a,b}$ avoids $v_a$, so $\Gamma_{a,b}$ is a smooth conic. As $b$ varies, these conics form a pencil on $S_a$. Thus the restriction of $\xi\circ\sigma$ to a general member of $\mathcal P$ is nonconstant, proving that $\phi$ is dominant. Since $\Gamma_{a,b}$ meets $\ell\cup C_4$ in finitely many points, we have $\overline{(S_a\cap S'_b)\setminus(\ell\cup C_4)}=\Gamma_{a,b}.$\\

\noindent The closure of a general fibre of $\phi$ is the smooth conic $\Gamma_{a,b}$. After resolving $\psi$ equivariantly, a general fibre is a smooth projective curve birational to this conic, hence isomorphic to $\pl$. The relative $G$-MMP gives a $G$-Mori fibre space over a surface, so $X$ is not $G$-solid.
\end{proof}

\begin{proof}[Proof of Proposition~\ref{prop:notgsol}]
Suppose that $\pic^G(X)\cong\dz$ and $G\not\cong\mathfrak A_5\times\dz_2$. In cases~\ref{joe1} and~\ref{joe2} of Theorem~\ref{theo:autQC}, Corollary~\ref{coro:finite-subgroups-autX} leaves the groups $G_{n,j}$ with $3\le j\le7$ and, in case~\ref{joe1}, the groups with $H\cong\mathfrak A_4$ or $\mathfrak S_4$. Corollary~\ref{coro:gmsubgroups} treats the former, and Proposition~\ref{prop:a4s4-not-solid} treats the latter.

In case~\ref{joe3}, the rational quotient in the proof of Lemma~\ref{lemm:gaz2z2} is $G$-equivariant and has rational general fibre. Resolving this map and running the relative $G$-MMP gives a $G$-Mori fibre space over a surface. Cases~\ref{joe4} and~\ref{joe5} follow from Lemma~\ref{lemm:Gis22or222}, and case~\ref{joe6} from Lemma~\ref{lemm:Gis2}.
\end{proof}

\begin{proof}[Proof of Theorem~\ref{theo:rigid}]
Proposition~\ref{prop:A5xZ2rigid} gives $G$-superrigidity when $G\cong\mathfrak A_5\times\dz_2$. Superrigidity implies rigidity, which in turn implies solidity. Proposition~\ref{prop:notgsol} excludes solidity for every other group under the rank-one hypothesis.
\end{proof}

\section{Linearisation}

The goal of this section is to prove Theorem \ref{theo:lin}. We will start by recalling the definition of the Amitsur groups and the obstruction to linearisation it provides, and apply it to finite groups $G\subset\aut(X)$ containing a subgroup isomorphic to $\dz_2^3$ to show that they are not projectively linearisable.

\subsection{Amitsur groups}

We briefly recall the facts about Amitsur groups that we will need; see \cite[Appendix~A]{blanc2023finite}. Let $Y$ be a smooth projective variety with an action of a finite group $\Gamma$. There is a homomorphism
$$
\delta:\pic(Y)^\Gamma\longrightarrow H^2(\Gamma,\dc^\times)
$$
whose kernel consists of classes of line bundles admitting a $\Gamma$-linearisation, that is, an equivariant lift of the $\Gamma$-action to their total space which is linear on fibres. The {\em Amitsur group} of $(Y,\Gamma)$ is defined by
$$
\operatorname{Am}(Y,\Gamma):=\operatorname{Im}(\delta).
$$
By \cite[Theorem~A.1]{blanc2023finite}, it is an invariant of the $\Gamma$-equivariant birational class of $Y$.

Recall that the canonical bundle has a natural $\Gamma$-linearisation. In particular, if $\pic^\Gamma(Y)=\dz[K_Y]$, then $\operatorname{Am}(Y,\Gamma)=0$. Moreover, for a biregular action of $\Gamma$ on $\mathbb P^n$, the equality $\operatorname{Am}(\mathbb P^n,\Gamma)=0$ means that $\mathcal O_{\mathbb P^n}(1)$ admits a $\Gamma$-linearisation. Equivalently, the action on $\mathbb P^n$ is induced by a linear representation of $\Gamma$.

\begin{prop}\label{prop:z23notlin}
If $G$ contains a subgroup isomorphic to $\dz_2^3$, then its action on $X$ is not projectively linearisable. In particular its action is not linearisable.
\end{prop}
\begin{proof}
Let $A\subset G$ be a subgroup isomorphic to $\dz_2^3$. As seen above, $A$ does not fix a point on $X$.

The kernel of the action of $\aut(X)$ on $\pic(X)$ is $\aut(Q,C_4)\subset\pgl_2(\dc)$, which does not contain a subgroup isomorphic to $\dz_2^3$. Hence, $A$ acts non-trivially on $\pic(X)$. Since $\pic(X)$ has rank two and $-K_X$ is primitive, it follows that

$$
\pic^A(X)=\dz[-K_X].
$$
The natural $A$-linearisation of the anticanonical bundle therefore gives

$$
\operatorname{Am}(X,A)=0.
$$
This vanishing also follows from
\cite[Lemma~3.8]{sharma2026amitsurgroupsimprimitivefano}.

Assume on the contrary that there exists a $G$-equivariant birational map $\chi:X\dashrightarrow\mathbb P^3$, where $\mathbb P^3$ is endowed with a biregular $G$-action. Restricting to $A$ and using birational invariance, we obtain

$$
\operatorname{Am}(\mathbb P^3,A)=\operatorname{Am}(X,A)=0.
$$
Thus, the $A$-action on $\mathbb P^3$ is induced by a four-dimensional linear representation of $A$. Since $A$ is abelian, this representation is a direct sum of one-dimensional representations. In particular, $A$ fixes a point on $\mathbb P^3$.

Applying \cite[Proposition~5.3]{reichstein2000essential} to the $A$-equivariant rational map $\chi^{-1}:\mathbb P^3\dashrightarrow X$, we conclude that $A$ fixes a point on $X$, a contradiction.
\end{proof}

\subsection{Fixed points and $G$-linearisability}

\begin{prop}\label{prop:fixptlin}
    Let $G$ be a finite subgroup of $\aut(X)$. If $G$ fixes a point in $X$, then $G$ is linearisable.
\end{prop}

To prove Proposition~\ref{prop:fixptlin}, we start by giving a geometric analogue of the rationality criterion over non-closed fields given in \cite[Section 5]{kuznetsov2024rationality}. We caution that the following proof is taken almost verbatim from the paper.

Let us first give some notation. Fix $G\subset\aut(X)$ finite such that $X$ is $G$-Fano, and suppose that $x\in X$ is a $G$-fixed point. Since $G$ exchanges $E$ and $E'$, any $G$-fixed point in $E\cup E'$ lies in $E\cap E'$. Then we have two cases: either $x\in E\cap E'$ or $x\notin E\cap E'$. When $x\in E\cap E'$ then there exist precisely two (anticanonical) lines $L,L'\subset X$ passing through $x$. Their union $L\cup L'$ is $G$-irreducible. On the other hand when $x\notin E\cap E'$, then no line on $X$ contains $x$. We define:
$$\Theta(x):=\begin{cases} L\cup L', \quad\text{when } x\in E\cap E',\\ \{x\}, \quad\text{when } x\notin E\cap E'.\end{cases}$$
and
$$m(x):=\begin{cases} 1, \quad\text{when } x\in E\cap E',\\ 2, \quad\text{when } x\notin E\cap E'.\end{cases}$$

\begin{lemm}
\label{lemma:Sarkisov-link}
    In the above setting, there exists a $G$-equivariant Sarkisov link
    \begin{equation}
    \label{eq:KP-link}
        \begin{tikzcd}
        & {\widetilde X} && {\widetilde X^+} \\
        X && {\bar X} && {X^+}
        \arrow[dashed, from=1-2, to=1-4, "\psi"]
        \arrow[from=1-2, to=2-1, "\sigma"']
        \arrow[from=1-2, to=2-3, "\phi"]
        \arrow[from=1-4, to=2-3, "\phi_+"']
        \arrow[from=1-4, to=2-5, "\sigma_+"]
        \end{tikzcd}
    \end{equation}
    where $\sigma$ is the blow-up of $X$ along $\Theta(x)$, $\phi$ and $\phi_+$ are small birational morphisms, $\psi$ is a flop and $\sigma_+$ is a $G$-elementary $K$-negative divisorial contraction, meaning that $\rk\,\pic(\widetilde X^+/X^+)^G=1$.

    When $x\notin E\cap E'$, then the threefold $X^+$ is a smooth del Pezzo threefold of degree 5, whereas when $x\in E\cap E'$ the threefold $X^+$ is a smooth divisor $V\subset \P^2\times\P^2$ of bidegree $(1,1)$, which is a del Pezzo threefold of degree 6. In this case, $V$ has a $G$-fixed point.
\end{lemm}

\begin{proof}
We first prove the existence of the link, and following this we identify $X^+$ and $\sigma_+$. Let $H=-K_X$, $m=m(x)$ and $F$ the $\sigma$-exceptional divisor. Since $\pic(X)^G=\mathbb Z[-K_X]$ and $F$ is $G$-irreducible, we have
$$\cl(\wt{X})_{\mathbb Q}^G=\mathbb Q H\oplus\mathbb Q F.$$
Both $H$ and $F$ are Cartier, so $\wt{X}$ is $G$-$\mathbb Q$-factorial and $\rk\,\pic(\wt{X})^G=2$. By \cite[Chapter 2, Theorems 2.2 and 3.4]{Isk}, we have that $H$ is very ample on $X$ and moreover the anticanonical image of $X$ is cut out by quadrics. Therefore if $\Theta(x)$ is the union of the two lines, we have that:
$$\Bs|H\otimes\mathcal{I}_{\Theta(x)/X}|=X\cap\langle\Theta(x)\rangle=\Theta(x),$$
which implies that $|H-F|$ is basepoint-free, whereas if $\Theta(x)=\{x\}$ then \cite[Lemma 5.6]{KP1} shows that $X\cap T_xX$ is the first-order neighbourhood of $x$, and hence that $|H-2F|$ is basepoint-free. This implies in either case that the anticanonical linear system $|-K_{\wt{X}}|$ of the blowup is basepoint-free.

Let $\oldphi\colon \wt{X}\to \P^{14-m}$ be the corresponding morphism, where the dimension of the target is calculated in \cite[Lemma 5.4(i) and (iii)]{KP1}. We have the following intersection numbers on $\wt{X}$:
\begin{equation}
\label{eq:Xtilde-intersections}
    H^3=28,\quad H^2\cdot F=0,\quad H\cdot F^2 =2(m-2),\quad F^3=m-1.
\end{equation}
In \cite[Lemma 5.4(iii) and (i)]{KP1} it is shown that the linear system $\mathscr{M}:=|H-(m+1)F|$ is positive dimensional. Let us show that $\mathscr{M}$ is mobile. For this, write $\mathscr{M}=\mathscr{N}+P$, where $\mathscr{N}$ is a mobile $G$-invariant linear system and $P$ is a nonzero $G$-invariant effective divisor. Pushing-forward to $X$ gives:
$$-K_X\sim\sigma_*(H-(m+1)F)\sim\sigma_*\mathscr{N}+\sigma_*P,$$
and since $\sigma_*\mathscr{N}$ is a positive-dimensional $G$-invariant linear system on $X$, and $\pic(X)^G=\mathbb Z[-K_X]$, we must have that $\sigma_*P=0$, and hence that $P$ is supported on $F$. Since $F$ is $G$-irreducible, we therefore have that $P=rF$ for some $r\geq1$.

The fact that $|H-mF|$ is basepoint-free and $\mathscr{N}$ is mobile, in conjunction with the intersection numbers (\ref{eq:Xtilde-intersections}), shows that:
$$0\leq (H-mF)\cdot \mathscr{N}^2=28-(m+r+1)\Bigl((m^2-3m+4)r+m(m-3)^2+4\Bigr).$$
When $m=1$, this simplifies to $30-2(r+3)^2\le -2$, whereas when $m=2$ this simplifies to $28-2(r+3)^2\le -4$, hence in both cases we obtain a contradiction. Hence $\mathscr{M}$ is mobile.

Now using (\ref{eq:Xtilde-intersections}) we have that $(H-mF)^3 = 22$ when $m=1$ and $(H-mF)^3 = 20$ when $m=2$, which implies that $\oldphi$ is generically finite. Let $\phi$ be the part of the Stein factorisation of $\oldphi$ with connected fibres, so that $\phi$ is birational. We claim that $\phi$ is small. Let $R\subset \wt{X}$ be a curve contracted by $\phi$. Then $R$ is contracted by $\oldphi$, so that $(H-mF).R=0$; nefness of $H$ then implies that $F.R\geq0$. If $F.R=0$, this would imply that $H.R=0$ so that $R$ is $\sigma$-exceptional, but then $F.R<0$ since $-F$ is $\sigma$-ample. Thus $F.R>0$, so that:
\begin{equation}
\label{eq:M}
    (H-(m+1)F).R<0,
\end{equation}
and so $R$ is contained in $\Bs\mathscr{M}$. Since $\mathscr{M}$ is mobile, this implies that $\phi$ is small. Therefore $\phi$ is a small birational contraction, and so the $G$-equivariant Sarkisov link (\ref{eq:KP-link}) exists by standard $G$-equivariant MMP theory, and is uniquely determined by $\Theta(x)$ (cf. \cite[Proposition 5.3]{KP1}).

Let us show that the contraction $\sigma_+$ of the other extremal ray of $\overline{\mathrm{NE}}(\wt{X}^+)^G$ is birational and divisorial. Let $\mathscr{M}^+$ denote the strict transform on $\wt{X}^+$ of the linear system $\mathscr{M}$; it is mobile since $\mathscr{M}$ is. Moreover, any $M^+\in\mathscr{M}^+$ is Cartier since terminal threefold flops preserve Cartierness, and by \eqref{eq:M} $M^+$ is positive on the flopped curves. If $M^+$ were not nef, it would therefore be negative on the other extremal ray of $\overline{\mathrm{NE}}(\widetilde X^+)^G$. Therefore $M^+.C<0$ for every $\sigma_+$-exceptional curve $C$, which implies that $\operatorname{Ex}(\sigma_+)$ is contained in $\Bs\mathscr{M}^+$. Since $\mathscr{M}^+$ is mobile, we have that  $\sigma_+$ is small. But this is impossible for a $K$-negative contraction from a terminal Gorenstein threefold by \cite[Theorem 0]{Benveniste1985}, hence $M^+$ is nef.

Let $R^+$ be the strict transform of the subvariety of $\wt{X}$ spanned by conics meeting $\Theta(x)$ if $m=1$, or by twisted cubics through $x$ if $m=2$. By \cite[Lemmas~2.11 and~2.12]{kuznetsov2024rationality}, $R^+$ is a divisor with two irreducible components satisfying
\begin{equation}
    \label{eq:R}
    R^+\sim H^+-(m+2)F^+\sim2M^++K_{\widetilde X^+},
\end{equation}
where $H^+$ (resp. $F^+$) is the strict transform on $\wt{X}^+$ be the strict transform of $H$ (resp. $F$). 
For the strict transform $C^+$ of a general curve in either family, the proof of \cite[Lemma~5.4]{kuznetsov2024rationality} gives $-K_{\widetilde X^+}. C^+=1$ and $F^+. C^+\geq 1$.
Hence $$M^+. C^+=1-F^+. C^+\leq 0,$$
which implies that $M^+.C=0$ for every $\sigma_+$-exceptional curve $C\subset\wt{X}^+$. Then $R^+. C=K_{\widetilde X^+}. C<0$, so $C\subset\operatorname{Supp}(R^+)$.
Therefore $\sigma_+$ is birational and divisorial, with exceptional locus $R^+$, and $\sigma_+^*A^+\sim M^+$ for some ample Cartier divisor $A^+$ on $X^+$.

Pushing-forward the relation \eqref{eq:R} to $X^+$ shows that $-K_{X^+}\sim 2A^+$, so that $X^+$ is a $G$-Fano threefold of index 2. From the intersection numbers \eqref{eq:Xtilde-intersections} and the fact that, using the projection formula, any triple intersection of Cartier divisors involving $-K_{\wt{X}^+}$ is not changed by flops, we get $(A^+)^3=7-m$. Hence $X^+$ is smooth by the classification in \cite[Theorem 1.7]{G-Fano}.

Since $\rk\,\cl(\wt{X}^+)=5-m$, we obtain $\rk\,\cl(X^+)=3-m.$ In the case where $m=1$, we have that $\rk\,\cl(X^+)=2$, which implies that $X^+$ is a divisor on $\P^2\times\P^2$ of bidegree $(1,1)$. Moreover, the threefold $\wt{X}$ has a single ordinary double point, so that $\wt{X}^+$ also has a single ordinary double point by \cite[Theorem 6.15]{kollar1998birational}. This point is necessarily fixed by $G$, which implies that $X^+$ has a $G$-fixed point.
\end{proof}

\begin{lemm}
\label{lemma:(1,1)-divisors}
Let $X$ be a smooth divisor on $\P^2\times\P^2$ of bidegree $(1,1)$, and let $G\subset\aut(X)$ be a finite subgroup. If $X$ has a $G$-fixed point, then the $G$-action is linearisable.
\end{lemm}

\begin{proof}
    Up to a change of coordinates we may assume that $X$ is given by the equation $x_0y_0+x_1y_1+x_2y_2=0$. Since $\aut(X)$ acts transitively on $X$, after conjugating $G$ we may assume that the point $p=\big([1:0:0],[0:0:1]\big)\in X$ is fixed by $G$, where $x_0,x_1,x_2$ and $y_0,y_1,y_2$ are coordinates on the two copies of $\P^2$. The $G$-action lifts to an action on $\P^2\times\P^2$, and every element of $G$ is of the form $g_A$ or $g_A\circ\tau$ for some $A\in \mathrm{GL}_3(\C)$, where 
    \begin{align*}
        g_A\colon \P^2\times\P^2 &\to \P^2\times\P^2 \\
        \big([x],[y]\big)&\mapsto\big([Ax],[A^{-T}y]\big) 
    \end{align*}
    and
    $$\tau\colon \big([x_0:x_1:x_2],[y_0:y_1:y_2]\big)\mapsto \big([y_2:y_1:y_0],[x_2:x_1:x_0]\big).$$
    Moreover, $\tau$ fixes $p$, and $g_A$ fixes $p$ precisely when $(1,0,0)$ is an eigenvector of $A$ and $(0,0,1)$ is an eigenvector of $A^T$; that is, when $A$ is upper triangular. The affine open subset $U=\{x_2y_0\neq0\}\subset X$ is therefore $G$-invariant.

    This open subset $U$ is isomorphic to $V=\C^3$ via the map:
    $$
       \phi\colon \big([x],[y]\big)\mapsto \left( \frac{x_1}{x_2}, \,\frac{y_1}{y_0},\,
\frac{x_0}{x_2}-\frac{y_2}{y_0} \right)
    $$
    with inverse
    $$\phi^{-1}\colon (u,v,z)\mapsto \big([z-uv : 2u: 2],[2:2v:-z-uv]\big).$$
    Let us compute the action of $G$ on $V$. Elements of the form $g_A$ act by:
$$(u,v,z)\mapsto \left(
\frac{du+e}{f},\,
\frac{av-b}{d},\,
\frac{ad\,z+bd\,u+ae\,v+2cd-be}{df}
\right)$$
for an invertible upper triangular matrix:
$$A=\begin{pmatrix} a &b &c \\ 0 & d & e \\ 0 & 0 &f \end{pmatrix},$$
whereas elements $g_A\circ\tau$ act by:
$$
(u,v,z)\mapsto
\left(
\frac{dv+e}{f},\,
\frac{au-b}{d},\,
\frac{-ad\,z+ae\,u+bd\,v+2cd-be}{df}
\right).
$$

Therefore $G$ acts on $V$ via affine transformations. Therefore the point 
$$\eta=\frac{1}{|G|}\sum_{g\in G}g(0)\in V$$
is fixed by $G$, so that for each $g\in G$ the function:
\begin{align*}
    \rho(g)\colon V&\to V, \\
    x&\mapsto g(x+\eta) - \eta
\end{align*}
is a linear isomorphism of $V$, and so the function $g\mapsto \rho(g)$ expresses $V$ as a three-dimensional $G$-representation. Considering $V$ as such, the affine varieties $U$ and $V$ are $G$-equivariantly isomorphic, so that $X$ is $G$-birational to $\P(V\oplus\C)$, where $G$ acts on the extra summand $\C$ trivially.
\end{proof}

\begin{proof}[Proof of Proposition \ref{prop:fixptlin}]
Let $x\in X$ be a $G$-fixed point. If $\pic^G(X)\not\cong\dz$, then the morphism $\pi\colon X\to Q$ is $G$-equivariant, and $\pi(x)$ is a $G$-fixed point of $Q$. We conclude by stereographic projection. Assume $\pic^G(X)\cong\dz$ and suppose that $x\notin E\cap E'$. Then by Lemma \ref{lemma:Sarkisov-link} $X$ is $G$-birational to the smooth del Pezzo threefold $V_5$. Since $G$ has a fixed point on $X$, we have that $G\ncong\mathfrak{A}_5$, so that $V_5$, and hence $X$, is $G$-birational to $\P(V)$ for some four-dimensional $G$-representation $V$ by \cite[Example 7.2]{cheltsov2025equivariant}.

On the other hand if $x\in E\cap E'$ then by Lemma \ref{lemma:Sarkisov-link}, $X$ is $G$-birational to a smooth divisor $V\subset\P^2\times\P^2$ of bidegree $(1,1)$, which has a $G$-fixed point. Then the $G$-action on $V$, and hence on $X$, is linearisable by Lemma \ref{lemma:(1,1)-divisors}.
\end{proof}

\begin{coro}\label{coro:a4s4-linearisation}
Let $G$ be as in Lemma~\ref{lemm:a4s4-actions}. Then the action of $G$ on $X$ is linearisable if and only if $G=\mathfrak S^{\mathrm{tw}}$.
\end{coro}

\begin{proof}
The group $\mathfrak S^{\mathrm{tw}}$ fixes $p_+$ and $p_-$, so Proposition~\ref{prop:fixptlin} applies. Each of the two direct products in Lemma~\ref{lemm:a4s4-actions} contains $\dz_2^2\times\langle\tau\rangle\cong\dz_2^3$, where $\dz_2^2$ is the normal Klein four subgroup of $\mathfrak A$. Proposition~\ref{prop:z23notlin} excludes projective linearisability in both cases.
\end{proof}

\subsection{Linearisation of the groups $G_{n,j}$}
\label{subsec:gm-linearization}

For odd $n\ge3$, the group $G_{n,5}$ has no fixed point on $X$ by Corollary~\ref{coro:fixed-points-autX}. We shall nevertheless linearise its action by passing to another quadric. The same model treats the other families $G_{n,j}$ that occur in Theorem~\ref{theo:lin}.

Throughout this subsection, $Q$ is in case~\ref{joe1} or~\ref{joe2} of Theorem~\ref{theo:autQC}, and $G=G_{n,j}$ up to conjugacy, with $j\in\{3,4,5,6,7\}$. We prove linearisability for $j\in\{3,4,6,7\}$ and for $j=5$ with $n$ odd.

We first construct a common equivariant birational model. Write $\lambda=1-4s^2$, with $\lambda=0$ in case~\ref{joe1}. Then
$$Q=\{s^2x_0x_4-x_1x_3+(1-s^2)x_2^2=0\}, \qquad s\notin\{0,\pm1\}.$$
Choose the sign of $s$ so that the birational involution $\sigma=\pi\circ\tau\circ\pi^{-1}$ of $Q$ is
$$\sigma=[f_4:sf_3:s^2f_2:sf_1:f_0],$$
as in \cite[Example~1.3]{malbon2024automorphisms}; in case~\ref{joe1}, take $s=\frac12$. For $\eta\in\dc^*$, denote by $t_\eta\in\aut(Q,C_4)$ the lift of $\left(\begin{smallmatrix}1&0\\0&\eta\end{smallmatrix}\right)$, so that $R_n=t_{\zeta_n}$. Put
$$\mathcal N=\langle t_\eta\ (\eta\in\dc^*),B,\tau\rangle\cong(\Gm\rtimes\dz_2)\times\dz_2.$$
Thus $G\subset\mathcal N$, and $\mathcal N=\aut(X)$ in case~\ref{joe2}.

\begin{lemm}\label{lemm:gm-quadric-model}
There is an $\mathcal N$-equivariant birational map
$$X\dashrightarrow Q_0=\{xy+zw+h^2=0\}\subset\mathbb P^4_{[x:y:z:w:h]},$$
where the action on $Q_0$ is given by
$$\begin{aligned} t_\eta:[x:y:z:w:h] &\longmapsto[x:y:\eta z:\eta^{-1}w:h],\\ B:[x:y:z:w:h] &\longmapsto[y:x:w:z:h],\\ \tau:[x:y:z:w:h] &\longmapsto[-x:-y:z:w:h]. \end{aligned}$$
These transformations lift to a linear representation on $\dc^5$.
\end{lemm}

\begin{proof}
Consider $\varphi_+=[x_1:x_2]$ and $\varphi_-=[x_3:x_2]$. For each sign, compose
$$(\varphi_\pm,\varphi_\pm\circ\tau):X\dashrightarrow\mathbb P^1\times\mathbb P^1$$
with projection of the Segre quadric from $([0:1],[0:1])$. This defines a rational map $\psi:X\dashrightarrow\mathbb P^2_p\times\mathbb P^2_q$. It is $\mathcal N$-equivariant: the torus preserves the projection centres, $\tau$ exchanges the rulings, and $B$ exchanges the two quadrics.

Writing
$$a=\frac{x_1}{x_2},\quad b=\frac{x_3}{x_2},\quad c=\frac{f_3}{sf_2},\quad d=\frac{f_1}{sf_2},$$
we have
$$p=[a+c:c-a:ac],\qquad q=[b+d:d-b:bd].$$
Substitution in the equation of $Q$ shows that the image is
$$F_s=\{p_2q_2=\alpha p_0q_0+\beta p_1q_1\}\subset\mathbb P^2\times\mathbb P^2,\qquad \alpha=\frac{1-s}{2},\quad\beta=\frac{1+s}{2}.$$
Since $\alpha\beta\ne0$, this is a smooth divisor of bidegree $(1,1)$. The map $\psi$ is birational: the two projections of Segre quadrics are birational, and $a,b,c,d$ recover a general point of $Q$. Indeed, on $x_2=1$, the functions $a,b$ determine $x_1,x_3,x_0x_4$, and $c,d$ then determine $x_0,x_4$.

On $F_s$, the torus element $t_\eta$ scales $p_2$ by $\eta^{-1}$ and $q_2$ by $\eta$; $\tau$ changes the signs of $p_1,q_1$; and $B$ exchanges the factors. The involution $R_2\tau$ is central in $\mathcal N$, so its fixed locus is $\mathcal N$-invariant. Its unique conic component is the smooth conic
$$\Gamma=F_s\cap\{p_0=q_0=0\},$$
which is therefore $\mathcal N$-invariant.

Under the Segre embedding $F_s\subset\mathbb P^7$, project from the plane $\langle\Gamma\rangle$. The resulting equivariant map is given by
$$[W:U:V:C:D]=[p_0q_0:p_1q_0:p_0q_1:p_2q_0:p_0q_2].$$
Its image is the smooth quadric
$$\mathcal Q=\{CD=\alpha W^2+\beta UV\}\subset\mathbb P^4,$$
and its inverse is
$$[W:U:V:C:D]\dashrightarrow\bigl([W:U:C],[W:V:D]\bigr).$$
Finally, setting
$$[x:y:z:w:h]=[\sqrt{\beta}\,U:\sqrt{\beta}\,V:iD:iC:\sqrt{\alpha}\,W]$$
identifies $\mathcal Q$ with $Q_0$ and gives the stated action. The corresponding matrices satisfy the group relations on $\dc^5$, so this action lifts to a linear representation.
\end{proof}

\begin{prop}\label{prop:gm-linearization}
Suppose that $G$ is conjugate to one of
$$G_{n,3},\qquad G_{n,4},\qquad G_{n,6},\qquad G_{n,7}\quad(n\geq1),$$
or to $G_{n,5}$ with $n$ odd. Then $X$ is $G$-birational to $\mathbb P(V)$ for a four-dimensional linear representation $V$ of $G$. In particular, $G$ is linearisable.
\end{prop}

\begin{proof}
We work on the model $Q_0$ of Lemma~\ref{lemm:gm-quadric-model}. Recall that projection from a $G$-fixed point $p\in Q_0$ is a $G$-equivariant birational map
$$Q_0\dashrightarrow\mathbb P(\dc^5/\ell_p),$$
where $\ell_p\subset\dc^5$ is the invariant line representing $p$.

For $G=G_{n,3}$ or $G_{n,4}$, the point $[1:0:0:0:0]$ is fixed. For $G=G_{n,6}$, the point $[1:-1:0:0:1]$ is fixed. Projection therefore proves the assertion in these cases.

For the remaining cases, let $\zeta$ be a primitive $N$-th root of unity and write
$$\rho_{a,b}:[x:y:z:w:h]\longmapsto[\zeta^a x:\zeta^{-a}y:\zeta^b z:\zeta^{-b}w:h].$$
For every $k\in\mathbb Z$, consider the birational self-maps
$$\begin{aligned} \Phi_k:[x:y:z:w:h] &\dashrightarrow[x(z/w)^k:y(w/z)^k:z:w:h],\\ \Psi_k:[x:y:z:w:h] &\dashrightarrow[x:y:z(x/y)^k:w(y/x)^k:h]. \end{aligned}$$
They preserve $Q_0$, commute with $B$, and have inverses $\Phi_{-k}$ and $\Psi_{-k}$. Moreover,
$$\Phi_k\rho_{a,b}\Phi_k^{-1}=\rho_{a+2kb,b},\qquad \Psi_k\rho_{a,b}\Psi_k^{-1}=\rho_{a,b+2ka}.$$
Suppose first that $G=G_{n,5}$ with $n$ odd. Then $\rho=R_n\tau$ has order $2n$ and $G=\langle\rho,B\rangle$. Relative to $\zeta=\zeta_{2n}$, we have $\rho=\rho_{n,2}$. Choose $\varepsilon\in\{1,-1\}$ such that $\varepsilon\equiv n\pmod4$. Conjugating first by $\Phi_{(\varepsilon-n)/4}$ and then by $\Psi_{-\varepsilon}$ changes the weights as follows:
$$(n,2)\longmapsto(\varepsilon,2)\longmapsto(\varepsilon,0).$$
The resulting action fixes $[0:0:1:1:i]\in Q_0$. Projection from this point proves linearisability.

Finally, let $G=G_{n,7}=\langle R_{2n}\tau,B\rangle$. Relative to $\zeta=\zeta_{2n}$, its rotation $\rho=R_{2n}\tau$ has weights $(n,1)$. Conjugation by $\Phi_{-\lfloor n/2\rfloor}$ changes these to $(0,1)$ if $n$ is even and $(1,1)$ if $n$ is odd. In the even case, the resulting action fixes $[1:1:0:0:i]$, and projection again gives a linearisation.

In the odd case, use the invariant affine chart $h=1$ and put
$$u=z-ix,\qquad v=w-iy,\qquad r=z+ix,\qquad q=w+iy.$$
The equation of $Q_0$ becomes $rv+uq+2=0$. On the open set $uv\ne0$, set
$$t=\frac{rv+1}{uv}.$$
Since
$$r=ut-\frac1v,\qquad q=-vt-\frac1u,$$
the functions $(u,v,t)$ are birational coordinates on $Q_0$. In these coordinates,
$$\rho:(u,v,t)\longmapsto(\zeta u,\zeta^{-1}v,t),\qquad B:(u,v,t)\longmapsto(v,u,-t).$$
This is a linear action on $\mathbb A^3$. Its standard projective compactification is the projectivisation of a four-dimensional representation, as required.
\end{proof}

\begin{proof}[Proof of Theorem~\ref{theo:lin}]
Linearisability implies projective linearisability by definition. By Proposition~\ref{prop:z23notlin}, projective linearisability implies that $G$ contains no subgroup isomorphic to $\dz_2^3$. It remains to prove linearisability under this last assumption.

Suppose first that $Q$ is in case~\ref{joe1} or~\ref{joe2} of Theorem~\ref{theo:autQC}. The groups $\aq\times\dz_2$, $\sq\times\dz_2$, $\ac\times\dz_2$, and $G_{n,5}$ with $n$ even all contain a subgroup isomorphic to $\dz_2^3$. Thus, by Corollary~\ref{coro:finite-subgroups-autX} and the assumption $\pic(X)^G\cong\dz$, the group $G$ is conjugate to one of $G_{n,3},G_{n,4},G_{n,6},G_{n,7}$, to $G_{n,5}$ with $n$ odd, or to $\mathfrak S^{\mathrm{tw}}$. The first five possibilities are linearisable by Proposition~\ref{prop:gm-linearization}. The case of $\mathfrak S^{\mathrm{tw}}$ follows from Corollary~\ref{coro:a4s4-linearisation}.

In cases~\ref{joe3}, \ref{joe5}, and~\ref{joe6}, every finite subgroup fixes a point on $X$ by Corollary~\ref{coro:fixed-points-autX}. In case~\ref{joe4}, our assumption implies that $G$ is a proper subgroup of $\aut(X)\cong\dz_2^3$, so the same corollary again gives a fixed point. Proposition~\ref{prop:fixptlin} therefore proves linearisability in all remaining cases.
\end{proof}

\bibliography{bib}
\bibliographystyle{alpha}
\end{document}

%% file: antoinescommands.tex
\newcommand{\dz}{\mathbb Z}
\newcommand{\dc}{\mathbb C}

\newcommand{\ac}{\mathfrak A_5}
\newcommand{\aut}{\mathrm{Aut}}
\newcommand{\mls}{\mathcal M}
\newcommand{\q}{\mathbb Q}
\newcommand{\cl}{\mathrm{Cl}}
\newcommand{\rk}{\mathrm{rk}}

\renewcommand{\phi}{\varphi}

\renewcommand{\P}{\mathbb{P}}

\newcommand{\sq}{\mathfrak S_4}
\newcommand{\pgl}{\mathrm{PGL}}

\newcommand{\pp}{\mathbb P^2}

\newcommand{\pl}{\mathbb P^1}

\newcommand{\pic}{\mathrm{Pic}}

\newcommand{\aq}{\mathfrak A_4}

\newcommand{\dih}{\mathrm{D}}

\newcommand{\p}{\mathbb P}

\newcommand{\mult}{\mathrm{mult}}
\renewcommand{\epsilon}{\varepsilon}
\newcommand{\dg}{\mathbb G}

\renewcommand{\tilde}{\widetilde}

\newtheorem{theo}[equation]{Theorem}
\newtheorem{prop}[equation]{Proposition}
\newtheorem{lemm}[equation]{Lemma}
\newtheorem{coro}[equation]{Corollary}
\newtheorem{rema}[equation]{Remark}